\documentclass[fleqno]{article}

\usepackage{amssymb}
\usepackage{amsthm}
\usepackage{amsmath}

\usepackage{color}

\def \tilde{\widetilde}

\newcommand{\ZZ}{{\mathbb Z}}

\newcommand{\Res}{{\rm Res}}

\newtheorem{theorem}{Theorem}%%%[section]
\newtheorem*{theorem*}{Theorem}
\newtheorem{lemma}{Lemma}
\newtheorem{proposition}{Proposition}
\newtheorem*{proposition*}{Proposition}

\theoremstyle{definition}
\newtheorem*{definition*}{Definition}
\newtheorem{definition}{Definition}
\newtheorem*{comments*}{Comments}

\newtheorem*{example*}{Example}

\theoremstyle{remark}
\newtheorem{remark}{Remark}
\newtheorem*{remarks*}{Remarks}

\makeatletter
\def\@maketitle{\newpage
 \null
 \vskip 2em
 \begin{center}%
  {\Large\bf \@title \par}%
  \vskip 1.5em
  {\normalsize
   \lineskip .5em
   \begin{tabular}[t]{c}\@author
   \end{tabular}\par}%
  \vskip 2em

 \end{center}%
 \par
 \vskip 2.5em}
\makeatother

\begin{document}
\title{An approach to Quantum Conformal Algebra}

\author{Carina Boyallian and Vanesa Meinardi\thanks{%
     Ciem - FAMAF, Universidad Nacional de C\'ordoba - (5000), C\'ordoba, Argentina- CIT-I.A.P Ciencias Humanas, Universidad Nacional de Villa Mar\'ia-(5900), C\'ordoba, Argentina
\newline $<$cboyallian@unc.edu.ar - vmeinardi@unvm.edu.ar$>$.}}

 \maketitle

\begin{abstract} We aim  to explore if inside  a quantum vertex algebras, we can find the right notion of a quantum conformal  algebra.

\end{abstract}

\section{ Introduction}

Since the pioneering papers \cite{BPZ,Bo1}, there has been a great
deal of work towards understanding of the algebraic structure
underlying the notion of the {o\-pe\-ra\-tor} product expansion
(OPE) of chiral fields of a conformal field theory. The singular
part of the OPE encodes the commutation relations of fields, which
leads to the notion of a  conformal algebra \cite{K}.

In \cite{BK}, they develop foundation of the theory of field algebras, which are a ``non-commutative version'' of a vertex algebra. Among other results they show that inside certain field algebras, more precisely strong field algebras ( where the $n$-product axiom holds)  we have a conformal algebra and a diferential algebra toghehter with certain compatibility equations, and conversely, having this two structures plus those equations we can recover a strong field algebra. One of these equations is the conformal analog of the Jacobi Identiy. They call a conformal algebra satisfying this equation {\it Leibnitz} conformal algebra.

A definition of a quantum vertex algebra, which is a deformation of a vertex algebra, was introduced by Etingof and Kazhdan in 1998,\cite{EK5}. Roughly speaking,  a quantum vertex algebra is a braided state-field correspondence
which {sa\-tis\-fies} associativity and braided locality axioms. Such braiding is a one-parameter braiding with coefficients in Laurent series.

Recently in \cite{DGK}, they developed a structure theory of quantum vertex algebras, parallel to that of vertex algebras.
In particular, they introduce braided n-products for a braided state-field correspondence and prove for quantum vertex algebras a version of the Borcherds
identity.

Following \cite{BK}, in this article, we try to determine the quantum analog of the notion of conformal algebra inside a quantum vertex algebra $V$. For this purpose, we introduced  new products parametrized by  Laurent polinomials $f$, and we showed that all this products are determined by those corresponding $f=1$ and $f=z^{-1}$.  The case $f=1$ coincides with the $\lambda$-product defining a conformal algebra(\cite{K},\cite{BK}). This allows us to deal with the coefficients of the braiding in $V$. An important remark is that $V$ together with the $\lambda$-product is no longer a Leibnitz conformal algebra, since due to the braiding, the analog of the Jacoby identity involves not only the products corresponding to $f=1$ (as in $\cite{BK}$), but those of $f= z^{-1}$. We translate to this language the hexagon axiom, quasi-associativity and associativity relations, and the braided skew-symmetry in a quantum vertex algebra, and all this allows us to give  an equivalent definition of quantum vertex algebra and present a candidate of a quantum conformal algebra.

The article is organized as follows. In Section 2 we review all the definitions and basic notion of field algebras and braided field algebras. In Section 3 we introduce the $(\lambda, f)$-product and prove some of its properties and we finish the section proving in Theorem 3 that shows that having a  strong braided field algebra is the same of having a conformal algebra, a differential algebra with unit with some compatibility equations.
In Section 4, we translate the hexagon axiom, quasi-associativity,\ and associativity relations, and the braided skew-symmetry in a quantum vertex algebra, we give  an equivalent definition of quantum vertex algebra and present a candidate of a quantum conformal algebra with an example of a Lie quantum conformal algebra.

\section{Preliminaries}

In this section review some basic definitions followig \cite{BK},\cite{DGK}. Throughout the paper all vector spaces, tensor products,etc are over a field $\mathbb{K}$ of characteristic zero, unless otherwise specified.

\subsection{Calculus of formal distribution}
Given a vector space $V$, we let $V[[z, z^{-1}]]$ be the space of  formal power series with coefficients in $V;$ they are called \textit{formal distributions}. A \textit{quantum field} over $V$ is a formal distribution $a(z) \in (\hbox{End}V)[[z, z^{-1}]]$ with coefficients in $\hbox{End}V,$ such that $a(z)v \in V((z))$ for every $v \in V.$ Hereafter $V((z))=V[[z]][z^{-1}]$ stands for the space of Laurent series with coefficients in $V.$

%Recall that the \textit{formal delta distribution} $\delta(z, w)$ is a formal distribution in $z$ and $w$ with coefficients in $\mathbb{K}$ defined as follows:
%\begin{equation}
%\delta(z, w)=\sum_{m \in \mathbb{Z}}z^{-m-1}w^m.
%\end{equation}
%It can be written as
%\begin{equation}
%\delta(z,w)=\iota_{z, w} \frac{1}{z-w}-\iota_{w,z}\frac{1}{z-w},
%\end{equation}

Throughout the article $\iota_{z,w}$ (resp $\iota_{w, z}$) denotes the geometric series expansion in the domain $|z|> |w|$ (resp $|w|>|z|$), namely we set for $n \in \mathbb{Z},$ $$\iota_{z,w}(z+w)^n=\sum_{l \in \mathbb{Z}_+}{n \choose l}z^{n-l}w^l$$
where $${n \choose l}=\frac{n(n-1)\cdots (n-l+1)}{l!}.$$

For an arbitrary formal distribution $a(z),$ we have
\begin{equation}
\hbox{Res}_z(a(z))=a_{-1},
\end{equation}
which is the coefficient of $z^{-1}.$
Denote by $\hbox{gl}f(V)$ the space of all $\hbox{End}V$-valued fields. We also need the \textit{Taylor's Formula} (cf. Proposition 2.4,[K1]), namely,

\begin{equation}\label{Taylor formula}
\iota_{z,w} a(z+w)=\sum _{j\in\ZZ_+}\frac{\partial_z^j}{j!}a(z) w^j=e^{w\partial_z}a(z).
\end{equation}

For each $n \in \mathbb{Z}$ one defines the $n$-th \textit{product} of fields $a(z)$ and $b(z)$ by the following formula:
\begin{equation}\label{n-product}
a(z)_{(n)}b(z)=\hbox{Res}_x(a(x)b(z)\iota_{x,z}(x-z)^n-b(z)a(x)\iota_{z,x}(x-z)^n).
\end{equation}
%XXXXXXXXXXXX HACE FALTA LO QUE SIGUE?
%Formula (\ref{n-product}) is equivalent to the following two formulas for $n \in \mathbb{Z}_{+}:$
%\begin{eqnarray*}
%  a(z)_{(n)}b(z) &=& \hbox{Res}_x[a(z), b(x)](x-z)^n ,\\
 % a(z)_{(-n-1)} b(z)&=& :\partial_z^n a(z)b(z):/n!.
%\end{eqnarray*}
%Here $:\,:$ stands for the \textit{normal ordered product} of field defined by
%\begin{equation}\label{normal ordered product}
%:a(z)b(z):=a(z)_{+}b(z)+b(z)a(z)_{-},
%\end{equation}
%XXXXXXXXXXXXXXXXXXXXXXXX
Denote by
$$a(z)_{+}=\sum_{j \leq -1} a_{(j)} z^{-j-1}, \qquad a(z)_{-}=\sum_{j \geq 0}a_{(j)} z^{-j-1}.$$

\
 \subsection{Conformal algebras and Field Algebras}
 In this subsection we recall the definition of a field algebra, conformal algebras and its properties following \cite{BK}

A \textit{state-field correspondence} on a pointed vector space $(V, |0\rangle)$ is a linear map $Y: V\otimes V \rightarrow V((z)), \, a\otimes b \rightarrow Y(z)(a\otimes b)$ satisfying
 \begin{enumerate}
 \item [(i)](vacuum axioms )$Y(z)( |0\rangle\otimes a)=a, \, Y(z)(a\otimes  |0\rangle)\in a+V[[z]]z;$
 \item [(ii)] (translation covariance) $TY(z)(a\otimes b)-Y(z)(a\otimes Tb)=\partial_zY(z)(a\otimes b)$,
 \item [(iii)]\hskip 3.5cm\, $Y(z)(Ta\otimes b)=\partial_z Y(z)(a\otimes b)$,
  \end{enumerate}
 where $T(a):=\partial_z(Y(z)(a\otimes |0\rangle))\mid_{z=0}=a_{(-2)}|0\rangle$, is called the translation operator.

    % \item [(iii)](translation covariance)\, $Y(z)(Ta\otimes b)=\partial_z Y(z)(a\otimes b).$

 Note that we will also denote by $Y$  the map $Y: V \rightarrow \hbox{End} V[[z, z^{-1}]], \, a \mapsto Y(a,z)=\sum_{k \in \mathbb{Z}} a_{(k)}z^{-k-1},$ such that $Y(a,z)b=Y(z)(a\otimes b).$

 Note that $Y(a,z)$ is a quantum  field, i.e $Y(a, z)b \in V((z))$ for any $b \in V.$

 The following results, proved in \cite{BK}, will be usefull in the sequel.

 \begin{proposition}(cf.  \cite{BK}, Prop.2.7). Given $Y: V\otimes V \rightarrow V((z))$ satisfaying conditions (i) and (ii) above, we have:
 \begin{enumerate}
 \item [(a)] $Y(z)(a\otimes  |0\rangle))=e^{zT}a;$
 \item [(b)] $e^{wT}Y(z)(1\otimes e^{-wT})=\iota_{z,w}Y(z+w).$

 If, moreover, $Y$ is a state-field correspondence, then
 \item[(c)] $Y(z)(e^{wT}\otimes 1)=\iota_{z,w}Y(z+w).$
 \end{enumerate}
 \end{proposition}

 Given a state field correspondence $Y,$ define
  \begin{equation}\label{Yop} Y^{op}(z)(u\otimes v)=e^{zT}Y(-z)(v\otimes u).\end{equation}
Then $Y^{op}$ is also a state-field correspondence, called the \textit{opposite} to $Y.$ (cf. [BK], Prop 2.8).

 Let $(V, |0\rangle)$ be a pointed vector space and let $Y$ be a state-field correspondence. Recall that $Y$ satisfies the $n$-th \textit{product axiom} if for all $a, b \in V$ and $n\in \mathbb{Z}$
 \begin{equation}\label{product axiom}
 Y(a_{(n)}b, z)=Y(a,z)_{(n)}Y(b,z).
 \end{equation}

 \,

 We say that $Y$ satisfies the \textit{associativity axiom} if for all $a, b, c \in  V$,  there exists $N\gg 0$ such that
 \begin{equation}\label{associativity axiom}
 \begin{aligned}
 (z-w)^NY(-w)&((Y(z)\otimes 1))(a\otimes b \otimes c)\\
 &=(z-w)^N\iota_{z,w}Y(z-w)(1\otimes Y(-w))(a\otimes b \otimes c).\end{aligned}
 \end{equation}

 \,

  Let $(V, |0\rangle)$ be a pointed vector space. As in \cite{BK}, a \textit{field algebra} $(V, |0\rangle, Y)$ is a state-field correspondence $Y$ for $(V, |0\rangle)$ satisfying the associativity axiom (\ref{associativity axiom}). A \textit{strong field algebra} $(V, |0\rangle, Y)$ is a state-field correspondence $Y$ satisfying the $n$-th product axiom (\ref{product axiom}). Note that, every strong field algebra is a field algebra.

  \

Let $(V, |0\rangle)$ be a pointed vector space and let $Y$  be a state-field corresponcence. For $a, b \in V, $ \cite{BK} defined the  $\lambda$-\textit{product} given by
\begin{equation}\label{lambda} a_{\lambda}b=\hbox{Res}_z e^{\lambda z}Y(z)(a\otimes b)=\sum_{n\geq 0}\frac{\lambda ^n}{n!} a_{(n)}b.
\end{equation}
and the $\cdot$-\textit{product} on $V$, which is denote as
\begin{equation}\label{dot} a\cdot b=\hbox{Res}_z z^{-1}Y(z)(a\otimes b)=a_{(-1)}b.\end{equation}
The vacuum axioms for $Y$ implies
\begin{equation}\label{vacum as unit}
|0\rangle  \cdot a=a =a \cdot |0\rangle,
\end{equation}
while the translation invariance axioms imply
\begin{equation}\label{t derivation  of the dot product}
T(a\cdot b)=T(a)\cdot b + a \cdot T(b),
\end{equation}
and
\begin{equation}\label{T derivation of lambda product}
T(a_{\lambda}b)=(Ta)_{\lambda}b+a_{\lambda}(Tb), \qquad (Ta)_{\lambda}b=-\lambda a_{\lambda}b
\end{equation}
for all $a, b \in V.$\
 Notice that from these equations we can derive that  $T(|0\rangle)=0$ and $|0\rangle_{\lambda}a=0=a_{\lambda}|0\rangle$ for $a \in V.$

Conversely, if we are given a linear operator $T,$ a $\lambda$-product and a $\cdot$-product
 on $(V,|0\rangle)$, satisfying the above properties (\ref{vacum as unit})-(\ref{T derivation of lambda product}), we can reconstruct the state-field correspondence $Y$ by the formulas
 \begin{equation}\label{reconstruction of state field correspondence}
 Y(a,z)_{+}b=(e^{zT}a)_{\cdot}b, \qquad Y(a, z)_{-}b=(a_{-\partial_z}b)(z^{-1}),
 \end{equation}
where $Y(a, z)=Y(a,z)_{+}+Y(a,z)_{-}.$

\,

 A $\mathbb{K}[T]$-module $V,$ equipped with a linear map $V \otimes V \rightarrow \mathbb{K}\otimes V, \, a\otimes b \rightarrow a_{\lambda}b,$ satisfying (\ref{T derivation of lambda product}) is called a $(\mathbb{K}[T])$-\textit{conformal algebra}. On the other hand with respect to the $\cdot$-product, $V$ is a  $(\mathbb{K}[T])$-\textit{differential algebra} (i.e an algebra with derivation $T$) with a unit $|0\rangle.$

\,

 Summarizing, (Cf. \cite{BK}, Lemma 4.1), we have that,
 giving a state-field correspondence on a pointed vector space $(V, |0\rangle)$ is equivalent to provide $V$ with a structure of a $\mathbb{K}$[T]-conformal algebra and a structure of a $\mathbb{K}[T]$-differential algebra with a unit $|0\rangle.$

 \,

Now,  recall the following results. Later on, we will prove some analogous result for the braided environment.

 \begin{lemma}(\cite{BK}, Lemma 4.2) Let $(V, |0\rangle)$ be a pointed vector space and let $Y$ be a state-field correspondence. Fix $a, b, c \in V.$ Then the collection of $n$-th product identities
 $
 Y(a_{(n)}b, z)c=(Y(a,z)_{(n)}Y(b,z))c$ (for $n\geq 0$)
 implies
 \begin{equation}\label{jacobi identity}
 (a_{\lambda}b)_{\lambda+\mu}=a_{\lambda}(b_{\mu}c)-b_{\mu}(a_{\lambda}c),
 \end{equation}

 \begin{equation}\label{compatibility lambda product and dot product}
 a_{\lambda}(b_{\cdot}c)=(a_{\lambda}b)_{\cdot}c+b_{\cdot}(a_{\lambda}c)+\int_{0}^{\lambda} (a_{\lambda}b)_{\mu}c\, d\mu.
 \end{equation}
 The $(-1)$-st product identity $Y(z)(a_{(-1)}b\otimes c)=(Y(z)_{(-1)}Y(z))(a\otimes b \otimes c)$ implies
 \begin{equation}
 (a{\cdot}b)_{\lambda}c=(e^{T\partial_{\lambda}}a)_{\cdot}(b_{\lambda}c)+(e^{T\partial_{\lambda}}b)_{\cdot}(a_{\lambda}c)
 +\int_{0}^{\lambda}b_{\mu}(a_{\lambda-\mu}c)\, d\mu,
 \end{equation}
 \begin{equation}\label{associativity dot product}
   (a_{\cdot} b)_{\cdot} c-a_{\cdot}(b_{\cdot} c)=\left(\int_{0}^{T}d\lambda\, a\right){\cdot}(b_{\lambda}c)+\left(\int_{0}^{T}d\lambda \,b\right)_{\cdot}(a_{\lambda}c).
 \end{equation}
 \end{lemma}

 Identity (\ref{jacobi identity}) is called the (left) \textit{Jacobi identity}. A conformal algebra satisfying this identity for all $a, b, c \in V$ is called a (left) \textit{Leibnitz conformal algebra}. Equation (\ref{compatibility lambda product and dot product}) is known as the ``non-commutative'' Wick formula, while (\ref{associativity dot product}) is called the \textit{quasi-associativity} formula.

 %%XXXXXX EESTO LO USAMOS DESPUES?????XXXXXX

 %Notice that the right hand side of (\ref{associativity dot product}) is symmetric with respect to $a$ and $b,$ hence %(\ref{associativity dot product}) implies
 %\begin{equation}\label{left symmetric}
 %  a_{\cdot}(b_{\cdot}c)-b_{\cdot}(a_{\cdot}c)=(a_{\cdot}b-b_{\cdot}a)_{\cdot}c.
 %\end{equation}
 %An algebra satisfying (\ref{left symmetric}) for all $a, b, c \in V $ is called \textit{left-symmetric.} For such an algebra %$a_{\cdot}b-b_{\cdot}a$ is a Lie algebra bracket.

 %XXXXXXXXXXXXXXXXXXXXXXXXXXXXXXXXXXXXXX

Finally, we also recall the following result.

 \begin{theorem}(\cite{BK}, Theorem 4.4)
 Giving a strong field algebra structure on a pointed vector space $(V, |0\rangle)$ is the same as providing $V$ with a structure of Leibnitz $\mathbb{K}[T]$-conformal algebra and a structure of a $\mathbb{K}[T]$-differential algebra with a unit $|0\rangle,$ satisfying (\ref{compatibility lambda product and dot product})-(\ref{associativity dot product}).
 \end{theorem}

 Recall also the following result.

 \begin{theorem}(\cite{BK}, Theorem 6.3)
 A vertex algebra is the same as a field algebra $(V, |0\rangle, Y )$ for which
$Y = Y^{ op}$.
 \end{theorem}

 Therefore we may assume this as a definition of vertex algebra.

 \,

  %\subsection{Braided Algebras}
\subsection{Braided Field Algebras}
We will follow the notation and presentation introduced in \cite{DGK}.

Throughout the rest of the paper we shall work over the algebra $\mathbb{K}[[h]]$ of formal series in the variable $h$, and all the algebraic structures that we will consider are modules over $\mathbb{K}[[h]]$.

 A \textit{topologically free $\mathbb{K}[[h]]$-module} is isomorphic to $W[[h]]$ for some $\mathbb{K}$-vector space $W.$

Note that $W[[h]]\ncong W\otimes \mathbb{K}[[h]],$ unless $W$ is finite-dimensional over $\mathbb{K},$ and that the tensor product $U[[h]]\otimes_{\mathbb{K}[[h]]}W[[h]]$ of topologically free $\mathbb{K}[[h]]$-modules is not topologically free, unless one of $U$ and $W$ are finite dimensional. For any vector space $U$ and $W,$ the \textit{completed} tensor product by
\begin{equation}\label{completed tensor product}
  U[[h]]\hat{\otimes}_{\mathbb{K}[[h]]}W[[h]]:=(U\otimes W)[[h]]
\end{equation}
This is a completion in $h$-adic topology of $U[[h]]\otimes_{\mathbb{K}[[h]]}W[[h]].$

\,

Given a topologically free $\mathbb{K}[[h]]$-module $V,$ we let
\begin{equation}
  V_h((z))=\left\{ a(z)\in V[[z, z^{-1}]]\,\mid\, a(z)\in V((z))\, \hbox{mod }\, h^M\, \hbox{for every}\, M\in \mathbb{Z}_{\geq 0}\right\}.
\end{equation}
Namely, expanding $a(z)=\sum_{n \in \mathbb{Z}}a_{(n)}z^{-n-1},$ we ask that  $$\lim _{n \rightarrow +\infty}a_{(n)}=0$$
in $h$-adic topology.

%\section{Quantum fields}

\noindent Let $V$ be a topologically free $\mathbb{K}[[h]]$-module. Following \cite{DGK}, we call \linebreak $\hbox{End}_{\mathbb{K}[[h]]}V$-\textit{valued quantum field}  an $\hbox{End}_{\mathbb{K}[[h]]}V$-valued formal distribution $a(z)$ such that  $a(z)b \in V_h((z))$ for any $b \in V.$

\,

\noindent Later on, we will need the  following lemmas,  proved in \cite{DGK}(cf. Lemma 3.2 and 3.3).

\begin{lemma} Let $|0\rangle\in V$ and $T:V\rightarrow V$ be a $\mathbb{K}[[h]]$-linear map such that $T(|0\rangle)=0.$ Then for any $\hbox{End}_{\mathbb{K}[[h]]}V$-valued quantum field $a(z)$ such that $[T, a(z)]=\partial_z a(z)$ (translation covariance), we have
\begin{equation}
  a(z)|0\rangle=e^{zT}a=\sum_{k \geq 0}\frac{T^{k}a}{k!}\, z^k,
\end{equation}
where $a=\hbox{Res}_z z^{-1}a(z)|0\rangle.$
\end{lemma}

\

\begin{lemma}Let $T:V\rightarrow V$ be a $\mathbb{K}[[h]]$-linear map and let $a(z)$ be an $\hbox{End}_{\mathbb{K}[[h]]}V$-valued quantum field such that $[T, a(z)]=\partial_z a(z).$ We have
\begin{equation}
e^{wT}a(z)e^{-wT}=\iota_{z,w}a(z+w).
\end{equation}
\end{lemma}

\

%\subsection {Braided vertex algebra}

\

  Let $V$ be a topologically  free $\mathbb{K}[[h]]$-module, with a given non-zero vector $|0\rangle\in V$ ( vacuum vector) and a $\mathbb{K}[[h]]$-linear map $T:V\rightarrow V$ such that $T(|0\rangle)=0$ (translation operator). Again, following \cite{DGK},

  \begin{enumerate}\item[(a)] \textit{A topological state-field correspondence} on $V$ is a linear map
  \begin{equation}
  Y: V\hat{{\otimes}}V \rightarrow V_h((z)),
  \end{equation}
  satisfying

 \begin{enumerate}
 \item [(i)](vacuum axioms)
 $Y(z)(|0\rangle\otimes v) = v$ and

 $$Y(z)(v\otimes|0\rangle) \in  v+V[[z]]z,\hbox{ for all } z\in V;$$

 \item[(ii)](\textit{translation covariance})\begin{equation}\label{translation covariance}\partial_{z} Y(z)=TY(z)-Y(z)(1\otimes T) = Y(z)(T\otimes 1),
    \end{equation}

   \end{enumerate}
   \item[(b)] A \textit{braiding} on $V$ is a $\mathbb{K}[[h]]-$linear map
   \begin{equation}
     \mathcal{S}(z): V\hat{\otimes} V \rightarrow V\hat{\otimes} V \hat{\otimes}(K((z))[[h]])
 \end{equation}
 such that $\mathcal{S}=1+O(h).$
  \end{enumerate}

  \

  A  \textit{braided state-field correspondance}  is a quintuple $(V, |0\rangle, T, Y, \mathcal{S})$ where $Y$ is a topological state-field correspondance and $S$ is a braiding as above.

   \

 % Note that for the topological state- field correspondence, we  have the fact that $T$ is determined by the map $Y$ and the vacuum vector $|0\rangle$. Namely,
 % $$T(v)=v_{(-2)}|0\rangle=\hbox{Res}_z\, (z^{-2} Y(z)(v\otimes |0\rangle)).$$

  \

 We will use the following standard notation: given $n\geq 2$ and $i, j\in\{1,\cdots, n\},$ we let
  \begin{equation}
    S^{i,j}(z): V^{\widehat{\otimes} n} \rightarrow V^{\widehat{\otimes}n}\hat{\otimes}(\mathbb{K}((z))[[h]]),
  \end{equation}
  act in the $i$-th and $j$-th factors (in this order) of $V^{\widehat{\otimes}_ n},$ leaving the other factors unchanged.

  \

A  \textit{braided vertex algebra} is a quintuple $(V, |0\rangle, T, Y, \mathcal{S})$ where $Y$ is a topological state-field correspondance and $S$ is a braiding as above, satisfying the following $\mathcal{S}$-\textit{locality}: for every $a, b \in V$ and $M \in \mathbb{Z}_{\geq 0},$ there exists $N=N(a, b, M)\geq 0$ such that
   \begin{equation}\label{S-locality}\begin{aligned}(z-w)^NY(z)(1\otimes Y(w))&\mathcal{S}^{12}(z-w)(a\otimes b \otimes c)\\
&= (z-w)^N Y(w)(1\otimes Y(z))(b\otimes a \otimes c),
  \end{aligned}\end{equation}
  where  this equality holds $\hbox{mod } h^M$, for all $c \in V.$

  \

  Again, given a  topological state-field correspondence $Y$, set
  \begin{equation}\label{Yop} Y^{op}(z)(u\otimes v)=e^{zT}Y(-z)(v\otimes u).\end{equation}

It was shown in \cite{DGK}, Lemma 3.6, that in a braided vertex algebra $V$ we have
  \begin{equation}\label{skew-symmetry}
    Y(z)\mathcal{S}(z)(a\otimes b)=Y^{op}(z)(a\otimes b)
  \end{equation}
  for all $a, b \in V.$

 After the proof of this result,(cf. Remark 3.7, [DGK]) they point out that it is enough to have the  $\mathcal{S}$-locality (\ref{S-locality}) holding just  for $c=|0\rangle,$ to prove that  $YS=Y^{op}$ in a braided vertex algebra. We will use this remark later.

 \

 We recall at this point  two important Propositions for our sequel.

\
\begin{proposition}(\cite{EK5}, Prop. 1.1) Let $V$ be a braided vertex algebra. for every $a, b, c \in V $ and $M \in \mathbb{Z}_{\geq 0},$ there exists $N \geq 0$ such that
$$\iota_{z,w}((z+w)^NY(z+w)(1\otimes Y(w))\mathcal{S}^{23}(w)\mathcal{S}^{13}(z+w)(a\otimes b\otimes c))$$
\begin{equation}\label{hexagonal relation}
=(z+w)^NY(w)\mathcal{S}(w)(Y(z)\otimes 1)(a\otimes b \otimes c)\quad \hbox{mod}\, h^{M}.
\end{equation}

\

\end{proposition}
\begin{proposition}\label{modulo kerY}(\cite{DGK}, Proposition 3.9) Let $(V, |0\rangle, T, Y, \mathcal{S})$ be a braided vertex algebra. Extend $Y(z)$ to a map $V\hat{\otimes}V\hat{\otimes}(\mathbb{K}((z))[[h]])$ in the obvious way. Then, modulo $\hbox{Ker}\,Y(z),$ we have
\begin{enumerate}
\item [(a)] $\mathcal{S}(|0\rangle \otimes a)\equiv |0\rangle,$ and $\mathcal{S}(z)(|0\rangle \otimes a)\equiv |0\rangle \otimes a;$
\item[(b)] $[T\otimes 1, \,\mathcal{S}(z)]\equiv-\partial _z \mathcal{S}(z)$ (left shift condition);
\item[(c)] $[1\otimes T, \,\mathcal{S}(z)]\equiv \partial _z \mathcal{S}(z)$ (right shift condition);
\item [(d)] $[T\otimes 1 + 1\otimes T, \,\mathcal{S}(z)]\equiv 0$ ;
\item [(e)] $\mathcal{S}(z)\mathcal{S}^{21}(-z)=1$(unitary).

\
Moreover, we have the quantum Yang-Baxter equation:
\item [(f)]$\mathcal{S}^{12}(z_1-z_2)\mathcal{S}^{13}(z_1-z_3)\mathcal{S}^{23}(z_2-z_3)\equiv \mathcal{S}^{23}(z_2-z_3)\mathcal{S}^{13}(z_2-z_3)\mathcal{S}^{12}(z_1-z_2),$

modulo $\hbox{Ker}(Y(z_1)(1\otimes Y(z_2))(1^{\otimes 2}\otimes Y(z_3)(-\otimes -\otimes - \otimes |0\rangle)))$.
\end{enumerate}

\end{proposition}
\
\section{On the structure of braided state-field {corres\-pon\-den\-ce}}

As in \cite{BK}, we aim  to show that there are, inside  certain braided vertex algebras, a ``braided conformal algebra'' and a ``differential algebra'' satisfying some family of equation. Conversely, we will show that given such structures under some nice conditions, we can give some reconstruction theorem.

\,

  Let $(V, |0\rangle,T, Y, \mathcal{S})$ be a braided-state field correspondence. For $n \in \mathbb{Z},$ the \textit{quantum n-product} $Y(z)_{(n)}^{\mathcal{S}}Y(z)$ is defined as

  \begin{equation}\label{quantum n-product}\begin{aligned}
  (Y(z)_{(n)}^{\mathcal{S}}Y(z))(a\otimes b &\otimes c)=Res_{x}(\iota_{x,z}(x-z)^nY(x)(1\otimes Y(z))(a\otimes b \otimes c)\\
  -&\iota_{z,x}(x-z)^nY(z)(1\otimes Y(x))S^{12}(z-x)(b\otimes a \otimes c)).\\
 \end{aligned} \end{equation}
%and we define the $\mathcal{S}$-\textit{bracket}

\

%\
%\begin{eqnarray}\label{S bracket}
%[Y(x), Y(z)]_{\mathcal{S}}(a\otimes b \otimes %c)&=&\iota_{x,z}Y(x)(1\otimes Y(z))(a\otimes b \otimes c)\nonumber\\
%  &-&\iota_{z,x}Y(z)(1\otimes Y(x))S^{12}(z-x)(b\otimes a %\otimes c))
%\end{eqnarray}
 % \begin{lemma}(\cite{GK}, Lemma 5.4) The quantum n-product (\ref{quantum n-product}) satisfies the following properties:
%  \begin{enumerate}
%  \item[(a)] field condition $(a, b, c \in V)$
%  \begin{equation}\label{field condition}
%    (Y(z)_{(n)}^{\mathcal{S}}Y(z))(a \otimes b \otimes c) \in V_{h}((z))
%  \end{equation}
%  \item [(b)] vacuum condition $(a, b \in V)$
%  \begin{equation}\label{field condition}
%    (Y(z)_{(n)}^{\mathcal{S}}Y(z))(a \otimes b \otimes \boldsymbol{1}) \in V[[z]],
%  \end{equation}
%  \item [(c)] translation covariance
%  \begin{equation}\label{translation covariance for fields}
%  T\circ(Y(z)_{(n)}^{\mathcal{S}}Y(z))-(Y(z)_{(n)}^{\mathcal{S}}Y(z))\circ(1\otimes 1 \otimes T)=\partial_z(Y(z)_{(n)}^{\mathcal{S}}Y(z))
%  \end{equation}
%  \end{enumerate}
%  \end{lemma}

Note that this definition differs but is equivalent to the definition introduced in [DGK] when you ask $S$ to satisfies unitary relation $S_{21}(-x)S(x)=1$ (cf. proposition 3e) which holds in a braided vertex algebra, where we are going to work. Now, we have the following result.

  \begin{lemma}\label{T derivation of quantum n-product}Given $(V, |0\rangle, T, Y, \mathcal{S})$  a braided state-field correspondence satisfying the equations

  \begin{equation}\label{Tx1 y S} [T\otimes 1,\, \mathcal{S}(z)]=-\partial_z S(z),
  \end{equation}
  \begin{equation}\label{1xT y S}
  [1\otimes T, S(z) ]=\partial_z \mathcal{S}(z).
  \end{equation}
  The quantum n-product (\ref{quantum n-product}) satisfies the following equation
  \begin{equation}
  \partial_z(Y(a,z)^{\mathcal{S}}_nY(b,z))=(\partial_zY(a,z))_n^{\mathcal{S}}Y(b,z)+Y(a,z)_n^{\mathcal{S}}(\partial_zY(b,z)).
  \end{equation}
  \end{lemma}
  \begin{proof}
   Applying the definition of quantum $n$-product (\ref{quantum n-product}), using integration by parts and translation covariance (\ref{translation covariance}), the LHS becomes

  \begin{eqnarray}\label{parte 1}
%  % \nonumber % Remove numbering (before each equation)
  &\,&Res_{x}\iota_{x,z}\partial_z((x-z)^nY(x)(1\otimes Y(z)))(a\otimes b \otimes c)\nonumber\\
  &\,&-Res_x\iota_{z,x}\partial_z((x-z)^nY(z)(1\otimes Y(x))S^{12}(z-x))(b\otimes a \otimes c)\nonumber\\
   &\,&= Res_{x}\iota_{x,z}\partial_z(x-z)^{n}Y(x)(1\otimes Y(z))(a\otimes b \otimes c)\nonumber\\
  &\,& +Res_{x}\iota_{x,z}(x-z)^{n}Y(x)(1\otimes \partial_z)(1\otimes Y(z))(a\otimes b \otimes c)\nonumber\\
    &\,&-Res_x\iota_{z,x}\partial_z(x-z)^{n}Y(z)(1\otimes Y(x))S^{12}(z-x)(b\otimes a \otimes c)\nonumber\\
     &\,&-Res_x\iota_{z,x}(x-z)^n\partial_zY(z)(1\otimes Y(x))S^{12}(z-x)(b\otimes a \otimes c)\nonumber\\
      &\,&-Res_x\iota_{z,x}(x-z)^nY(z)(1\otimes Y(x))\partial_zS^{12}(z-x)(b\otimes a \otimes c)\nonumber\\
  &\,&= -Res_{x}\iota_{x,z}\partial_x(x-z)^{n}Y(x)(1\otimes Y(z))(a\otimes b \otimes c)\nonumber\\
  &\,& +Res_{x}\iota_{x,z}(x-z)^{n}Y(x)(1\otimes Y(z))(1\otimes T \otimes 1)(a\otimes b \otimes c)\nonumber\\
    &\,&+Res_x\iota_{z,x}\partial_x(x-z)^{n}Y(z)(1\otimes Y(x))S^{12}(z-x)(b\otimes a \otimes c)\nonumber\\
     &\,&-Res_x\iota_{z,x}(x-z)^nY(z)(T\otimes 1)(1\otimes Y(x))S^{12}(z-x)(b\otimes a \otimes c)\nonumber\\
      &\,&+Res_x\iota_{z,x}(x-z)^nY(z)(1\otimes Y(x))\partial_xS^{12}(z-x)(b\otimes a \otimes c)\nonumber\\
       &\,&= Res_{x}\iota_{x,z}(x-z)^{n}Y(x)(T\otimes 1)(1\otimes Y(z))(a\otimes b \otimes c)\nonumber\\
  &\,&+ Res_{x}\iota_{x,z}(x-z)^{n}Y(x)(1\otimes Y(z)(1\otimes T \otimes 1)(a\otimes b \otimes c)\nonumber\\
    &\,&-Res_x\iota_{z,x}(x-z)^{n}Y(z)(1\otimes Y(x))(1\otimes T\otimes 1)S^{12}(z-x)(b\otimes a \otimes c)\nonumber\\
     &\,&-Res_x\iota_{z,x}(x-z)^nY(z)(1\otimes Y(x))( T\otimes 1\otimes 1)S^{12}(z-x)(b\otimes a \otimes c).\nonumber\\
      \end{eqnarray}
      On the other hand using translation covariance, RHS becomes
      \begin{eqnarray}
      &\,&Res_x\iota_{x,z}(x-z)^nY(x)(T\otimes 1)(1\otimes Y(z))(a\otimes b \otimes c)\nonumber\\
      &\,&-Res_x \iota_{z,x}(x-z)^nY(z)(1\otimes Y(x))\mathcal{S}^{12}(z-x)(1\otimes T\otimes 1)(b\otimes a \otimes c)\nonumber\\
      &\,&+Res_x\iota_{x,z}(x-z)^nY(x)(1\otimes Y(z))(1\otimes T\otimes1)(a\otimes b \otimes c)\nonumber\\
       &\,&-Res_x \iota_{z,x}(x-z)^nY(z)(1\otimes Y(x))\mathcal{S}^{12}(z-x)(T \otimes 1\otimes 1)(b\otimes a \otimes c).\nonumber\\
      \end{eqnarray}

Due to equations (\ref{Tx1 y S}) and (\ref{1xT y S}) we get
 \begin{equation}\label{Tx1x1 y S}\begin{aligned}
  (T\otimes 1 \otimes 1)&\mathcal{S}^{12}(z-x)(b\otimes a \otimes c)=\\
  &\mathcal{S}^{12}(z-x)(T\otimes 1 \otimes 1)(b\otimes a \otimes c)+\partial_x\mathcal{S}^{12}(z-x)(b\otimes a \otimes c),
  \end{aligned}\end{equation}

  and

  \begin{equation}\label{1xTx1 y S}\begin{aligned}
  (1\otimes T \otimes 1)&\mathcal{S}^{12}(z-x)(b\otimes a \otimes c)=\\&\mathcal{S}^{12}(z-x)(1\otimes T \otimes 1)(b\otimes a \otimes c)-\partial_x\mathcal{S}^{12}(z-x)(b\otimes a \otimes c).
 \end{aligned} \end{equation}

  Applying  equations (\ref{Tx1x1 y S}) and (\ref{1xTx1 y S}) to RHS, we get
  \begin{eqnarray}\label{parte 2}
      &\,&Res_x\iota_{x,z}(x-z)^nY(x)(T\otimes 1)(1\otimes Y(z))(a\otimes b \otimes c)\nonumber\\
      &\,&-Res_x \iota_{z,x}(x-z)^nY(z)(1\otimes Y(x))(1\otimes T \otimes 1)\mathcal{S}^{12}(z-x)(b\otimes a \otimes c)\nonumber\\
       &\,&+Res_x \iota_{z,x}(x-z)^nY(z)(1\otimes Y(x))\partial_x\mathcal{S}^{12}(z-x)(b\otimes a \otimes c)\nonumber\\
      &\,&+Res_x\iota_{x,z}(x-z)^nY(x)(1\otimes Y(z))(1\otimes T\otimes1)(a\otimes b \otimes c)\nonumber\\
       &\,&-Res_x \iota_{z,x}(x-z)^nY(z)(1\otimes Y(x))(T \otimes 1\otimes 1)\mathcal{S}^{12}(z-x)(b\otimes a \otimes c)\nonumber\\
       &\,&-Res_x \iota_{z,x}(x-z)^nY(z)(1\otimes Y(x))\partial_x\mathcal{S}^{12}(z-x)(b\otimes a \otimes c)\nonumber\\
       &\,&=Res_x\iota_{x,z}(x-z)^nY(x)(T\otimes 1)(1\otimes Y(z))(a\otimes b \otimes c)\nonumber\\
      &\,&-Res_x \iota_{z,x}(x-z)^nY(z)(1\otimes Y(x))(1\otimes T \otimes 1)\mathcal{S}^{12}(z-x)(b\otimes a \otimes c)\nonumber\\
      &\,&+Res_x\iota_{x,z}(x-z)^nY(x)(1\otimes Y(z))(1\otimes T\otimes1)(a\otimes b \otimes c)\nonumber\\
       &\,&-Res_x \iota_{z,x}(x-z)^nY(z)(1\otimes Y(x))(T \otimes 1\otimes 1)\mathcal{S}^{12}(z-x)(b\otimes a \otimes c).\nonumber\\
       \end{eqnarray}
       Then equations (\ref{parte 1}) and (\ref{parte 2}) are equal, therefore  the claim follows.
       \end{proof}

       \

      \begin{remark} Recall that, as we quote in  Proposition \ref{modulo kerY}, it was shown by  \cite{DGK}  that in a braided vertex algebra, equations (\ref{Tx1 y S}) and (\ref{1xT y S}) hold mod $\hbox{Ker} Y$.  In \cite{EK5}, condition (\ref{Tx1 y S}) is asked as part of the definition of a braided vertex operator algebra. In this context,   asking (\ref{1xT y S}) is equivalent to ask $T$ to be a derivation of a braided vertex operator algebra. It is shown in \cite{Li}, that if in addition  we ask the undelying field algebra to be \textit{non-degenerate} (cf. definition 5.12, \cite{Li}), we have that (\ref{1xT y S}) holds in a braided vertex algebra where the associativity relation (\ref{associativity relation}) holds (cf. \cite{EK5}).\end{remark}

      \

  Let $(V, |0\rangle,T, Y, \mathcal{S})$ be a braided-state field correspondence.  $Y$ satisfies the \textit{quantum} $n$-th \textit{product identities} if for all $a, b, c \in V$ and $n \in \mathbb{Z}$
  \begin{equation}\label{axiom quantum n-product}
  Y(z)_{(n)}^{\mathcal{S}}Y(z)(a\otimes b \otimes c)=Y(z)(a^{\mathcal{S}}_{(n)}b\otimes c),
  \end{equation}
  where
  \begin{equation}
    a_{(n)}^{\mathcal{S}}b=\hbox{Res}(z^nY(z)S(z)(a\otimes b)).
  \end{equation}

\
%XXX PONEMOS ESTO XXXX
\

  $Y$ satisfies the \textit{ associativity relation} if for any $a, b, c \in  V$ and  $M \in \mathbb{Z}_{\geq 0}$ there exists $N \in \mathbb{Z}_{\geq 0}$ such that
 $$\iota_{z,w}(z+w)^NY(z+w)((1\otimes Y(w)))(a\otimes b \otimes c)$$
 \begin{equation}\label{associativity relation}
 =(z+w)^NY(w)(Y(z)\otimes 1)(a\otimes b \otimes c)\, \hbox{mod}\,h^M,
 \end{equation}

  Let $(V, |0\rangle)$ be a pointed vector space. We define  a \textit{ braided field algebra} $(V, |0\rangle, Y, T, \mathcal{S})$ is a braided state-field correspondence $Y$  satisfying the associativity relation (\ref{associativity relation}). We also introduce a \textit{strong braided field algebra} $(V, |0\rangle, Y, T, \mathcal{S})$ as a state-field correspondence $Y$ satisfying the quantum $n$-th product identities (\ref{axiom quantum n-product}). This are the braided versions of field algebra and strong field algebra introduced by \cite{BK}.

  \

Let $(V, |0\rangle, T, \mathcal{S})$ be a braided-state field correspondence. For $a, b \in V, f \in \mathbb{K}((z))[[h]],$ we define the  $(\lambda, f)$-\textit{product} by the
   \begin{equation}\label{(lambda, f)-product}
   a_{(\lambda, f)}b=\hbox{Res}_ze^{\lambda z}f(z)Y(z)(a\otimes b)=\sum_{ n \in Z_{\geq 0},\, \hbox{finite }}\sum_{i \in  \mathbb{Z}}f_i(h)\, a_{(n+i)}b\, \frac{\lambda^n}{n!} \in V\otimes \mathbb{K}[\lambda][[h]]
   \end{equation}
where $f(z)=\sum_{i\in \mathbb{Z}}f_i(h) z^{i},\, f_i(h) \in \mathbb{K}[[h]]. $ Note that $f_i(h)=0$ for $i<<0$.

\

\begin{remark}\begin{enumerate}
\item[(i)] If in addition we ask $V$ to have a structure of $\mathbb{K}((z))$-module structure, more precisely $z^k(a_{(n)}b)= a_{(n+k)}b$, this $(\lambda, f)$-product resembles the operations introduced in \cite{GKK}. Instead,  we are asking $V$ to have a braiding that involves some elements of  $\mathbb{K}((z))[[h]]$.
\item [(ii)] The product $a_{\lambda, f}$ coincides with $X_{\lambda, -\lambda-\partial}^{z_0, z_1}(a, b ; f(z_1-z_0))$ from Sec. 6.4 in the paper \cite{BDHK}.
\end{enumerate}
\end{remark}

%XXXXXXXXXXXXXXXXXX AGREGAR RELACION $(\lambda, z^{-1})$-product con $\cdot$-product XXXXXXXXXXXXXXXX
\

  % As in \cite{BK},  we denote
   %\begin{equation*}
  % a_{\cdot}b=\hbox{Res}_z z^{-1}Y(z)(a\otimes b)=a_{(-1)}b.
  % \end{equation*}

We have the following useful Lemma.

\begin{lemma} \label{relation between lambda product y lambda f product} Given $(V, |0\rangle, T, \mathcal{S})$ be a braided-state field correspondence, we have
 \begin{itemize}
       \item [(a)] $a_{(\lambda, z^mf)}b= \partial_{\lambda}^{m}a_{(\lambda,f)} b$ for $ \, m\geq 0$, and $f\in\mathbb{K}((z))[[h]]$. In particular, $a_{(\lambda, z^m)}b= \partial_{\lambda}^{m}a_{(\lambda,1)} b$ for $ \, m\geq 0$,
       \item [(b)] $a_{(\lambda, z^{-k})}b= (({\lambda+T})^{(k-1)}a)_{(\lambda, z^{-1})} b,$ for $\, k\geq 1$,
      % \item[(c)]$a_{(\lambda, z^{-1})}=a_{\cdot}b+\int_{0}^{\lambda}a_{\mu}b\, d\mu$.
     \end{itemize}
\end{lemma}
\begin{proof}
Let $f(z)=\sum_if_i(h)z^i,$ item $(a)$ follows from the definition of $(\lambda, f)$-product:
\begin{eqnarray*}
 a_{(\lambda, z^mf)}b&=&\hbox{Res}_ze^{\lambda z}z^m f(z)Y(z)(a\otimes b)\nonumber\\
 &=&\sum_if_i(h)\hbox{Res}_z\sum_{k\geq 0} \lambda^k/k!\sum_{j\in \mathbb{Z}} a_{(j)}b\,z^{-j-1+k+m+i}\nonumber\\
 &=& \sum_if_i(h)\sum_{k \geq 0}\lambda ^{k-m}/(k-m)!\, a_{(k+i)}b\nonumber\\
 %&=& \sum_{k \geq m}\lambda^{k-m}/(k-m)!a_{(k)}b\nonumber\\
 &=& \partial _{\lambda}^m a_{(\lambda, f)}b.
 \end{eqnarray*}
Applying definition of $(\lambda, f)$-product and using integration by parts and translation covariance we get item $(b)$, namely:
 \begin{eqnarray*}
 a_{(\lambda, z^{-k})}b&=&a_{(\lambda,\, (-\partial)^{(k-1)}_zz^{-1})}b\nonumber\\
 &=&\hbox{Res}_z e^{\lambda z}(-\partial)^{(k-1)}_zz^{-1}Y(z)(a\otimes b)\nonumber\\
 &=&\hbox{Res}_z z^{-1}\partial_z^{(k-1)} e^{\lambda z}Y(z)(a\otimes b)\\
 &=&\hbox{Res}_z z^{-1} e^{\lambda z}\sum_{r=0}^{k-1}\lambda^{(r)}Y(z)(T^{(k-1-r)}a\otimes b)\\
 &=&\hbox{Res}_z z^{-1} e^{\lambda z}Y(z)((\lambda+T)^{(k-1)}a\otimes b)\\
 &=& (({\lambda+T})^{(k-1)}a)_{(\lambda, z^{-1})} b.
 \end{eqnarray*}
%Item $(c)$ follows from definition and the fact that $z^{-1}e^{\lambda z}=z^{-1}+\int_{0}^{\lambda} e^{\mu}z\, d\mu.$
\end{proof}

    Note that if $f=1$ in (\ref{(lambda, f)-product}), we recover the $\lambda$-product introduced in (\ref{lambda}) for a state-field correspondence.
     We will denote $a_{(\lambda, 1)}=a_{\lambda}b$. Observe also that, due to the Lemma above, any $(\lambda,f)$-product can be written in terms of the $\lambda$-product and the  $(\lambda,z^{-1})$-product.

     \

     \
   The vacuum axioms for $Y$ imply that,
   \begin{equation}\label{identity}
   |0\rangle_{(\lambda,z^{-1})}a=a= a_{(\lambda,z^{-1})}|0\rangle,
   \end{equation}

   while the translation invariance axioms show that,
  % \begin{equation}\label{derivation of the dot product}
   %  T(a_{(\cdot,f)}b)=T(a)_{(\cdot,f)}b+a_{(\cdot,f)}T(b)
   %\end{equation}
   %and
   \begin{equation}\label{derivation of the lambda productprima}
     T(a_{(\lambda, f)}b)=T(a)_{(\lambda, f)}b+a_{(\lambda, f)}T(b)\end{equation} and
     \begin{equation}\label{derivation of the lambda product1prima}  T(a)_{(\lambda, f)}b=-\lambda a_{(\lambda, f)}b-a_{(\lambda, f')}b
   \end{equation}
   for all $a, b \in V$ and $f\in \mathbb{K}((z))$. Note that, when $f=1$ in (\ref{derivation of the lambda productprima}) and (\ref{derivation of the lambda product1prima}), we recover equation (\ref{T derivation of lambda product}).

   \,

   Conversely, if we are given a pointed topologically free $\mathbb{K}[[h]]$-module $(V,|0\rangle ),$ togheter with a $\mathbb{K}[[h]]$- linear map $T$, a braiding $\mathcal{S}$, a $(\lambda, 1)$-product and a $(\lambda, z^{-1})$-product on $V$ satisfying the properties (\ref{identity})-(\ref{derivation of the lambda product1prima}), we can reconstruct the braided state-field correspondence $Y$ by the formulas:
   \begin{equation}\label{definition state field correspondence}
    Y(a,z)_{+}b=(e^{zT}a)_{(\lambda, z^{-1})}b|_{\lambda=0}, \qquad Y(a,z)_{-}b=(a_{(-\partial_z,\,1)}b)(z^{-1})),
      \end{equation}
      where $Y(a,z)=Y(a,z)_{+}+Y(a,z)_{-}.$

      \

%XXXXXXXXXXREVISARXXXXXXXXXXXXX     A $\mathbb{K}[[T]]$-module $V,$ equipped with a linear map $V\otimes V \otimes\mathbb{K}((z))\rightarrow \mathbb{K}[[\lambda]]\otimes V,$ $a\otimes b\otimes f(z)\rightarrow a_{(\lambda,f)}b,$ satisfying (\ref{derivation of the lambda product}) is called a $\mathbb{K}[T]$-\textit{braided conformal algebra}.XXXXXXXXXXXXXX On the other hand with respect to the $\cdot$-product, $V$ is a $\mathbb{K}[T]$-differential algebra.

We will need the following Lemma.

 \begin{lemma}\label{Derivada de f en lambda prod}We have that $$a_{(\lambda, f^{(l)})}b=  ((-\lambda -T)^la)_{(\lambda, f)}b,$$
 for all $a$ and $b\in V$ and $l\geq 0$. Here and further $f^{(l)}(z)=\partial_z^lf(z)$.
 \end{lemma}

 \begin{proof} Straightforward using (\ref{derivation of the lambda product1prima}).\end{proof}

 For the following Proposition it will be useful to introduce the following notation:

   %\
  % XXX aca definimos el punto f producto?XXX

   \begin{equation}
   a_{(\cdot,f)}b:= a_{(\lambda, z^{-1}f)}b|_{\lambda=0}=\hbox{Res}_z z^{-1}f(z)Y(z)(a\otimes b)=\sum_{i \in  \mathbb{Z}}f_i(h)\, a_{(i-1)}b,
   \end{equation}
for $a, b \in V, f \in \mathbb{K}((z))[[h]], \, f(z)=\sum_{i\in \mathbb{Z}}f_i(h) z^{i}, f_i(h)\in \mathbb{K}[[h]]. $ Note that in the case $f=1$ we obtain the $\cdot-$product in [BK], ( cf. (\ref{dot})), namely

   \begin{equation*}
   a_{(\cdot,1)}b=a_{(\lambda,z^{-1})}b|_{\lambda=0}=a\cdot b,
   \end{equation*}

since it  is easy to show that

\begin{equation}\label{dot con integral}
a_{(\lambda, z^{-1})}b=a_{\cdot}b+\int_{0}^{\lambda}a_{\mu}b\, d\mu.\end{equation}
 \begin{remark} The $a_{\lambda, z^{-1}}$ from (\ref{dot con integral}) is the so called integral of $\lambda$-bracket, which appeared earlier in \cite{DK}.
 \end{remark}
Whith all this, we can state the following result.

   \begin{proposition}\label{equations of a braided conformal algebra}
  Let $(V, |0\rangle,T,Y, \mathcal{S})$ be a braided state field correspondence such $\mathcal{S}$-locality holds for $c=|0\rangle$ . Then the collection of the $n$-th quantum product identities (\ref{quantum n-product}) for $n\geq -1$ implies:

   \begin{equation}\label{jacoby identity (lambda,f) product}
   (a_{-\alpha-T}b)_{\alpha+\beta}c= -b_{\alpha}(a_{(\beta) }c)
      +\sum_{i=1}^{r}\sum_{l\geq 0}(-1)^la^{i}_{(\beta,(f_{i}(z))^{(l)})}(b^{i}_{(\alpha, x^l)}c),
   \end{equation}

\

   \begin{eqnarray}\label{relacion lambda y punto producto}
     (a_{-\lambda}b)_{\cdot}c &=&- b_{(\lambda-T)}(a_{\cdot}c)+\sum_{i=1}^{r}\sum_{l\geq 0}(-1)^l a^{i}_{(\cdot,(f_{i}(z))^{(l)})}(b^{i}_{(\lambda-T, x^l)}c)\nonumber\\
     &+&\int_{0}^{T-\lambda}(a_{-\lambda}b)_{\mu}c\,d\mu,
     \end{eqnarray}

      \begin{eqnarray}\label{relation dot product and lambda product}
     (a_{\cdot}b)_{\lambda}c &=& (e^{T\partial_{\lambda}}b)_{\cdot}(a_{\lambda}c) -\int_{0}^{-T}(a_{-\mu-T}b)_{\lambda}c\, d\mu
     + \sum_{i=1}^{r}[(e^{T\partial_{\lambda}}a^i)_{(\cdot, f_i(z))}(b^{i}_{\lambda}c)\nonumber\\
     &-&\sum_{l\geq 0}\int_{0}^{\lambda}a^i_{(\mu, (f_i(z))^{(l)})}(b^i_{(\lambda-\mu, x^l)}c\,)d\mu],
     \end{eqnarray}

  \begin{eqnarray}\label{asociatividad del producto punto}
   (a\cdot b)_{\cdot}c&=&b_{\cdot}(a_{\cdot}c)+\hbox{Res}_z\left(\int_{0}^{T}d_{\lambda}b\right)_{\cdot}(b_{\lambda}c)-\int_{0}^{-T}(a_{\mu-T}b)_{\cdot}c\,d\mu\nonumber\\
    &+&\sum_{i=1}^{r}\sum_{l\geq 0}\left(\int_{0}^{T}d_{\lambda}a^{i}\right)_{\cdot}(b^{i}_{\lambda,\, D_lf_i(z)}c)\nonumber\\
    &+&\sum_{i=1}^{r}\sum_{m,\,l\geq 0}(-1)^l a^{i}_{(\cdot, z^{m+1})}(b^{i}_{(.,\, D_{l}(f_i(z))z^{-m})}c),
   \end{eqnarray}
   where $D_l=z^l\partial_z^{(l)}$  and $\mathcal{S}(z)(a\otimes b)=\sum_{i=0}^{r} f_i(z)a^i\otimes b^i.$
   \end{proposition}
   \begin{proof}

Recall that the fact that the $\mathcal{S}$-locality holds for $c=|0\rangle$, implies that $Y(z)S(z)=Y^{op}(z)$. Applying definitions of $\lambda$-product and the definition of $Y^{op},$ due to Lemma \ref{skew-symmetry} we get
   \begin{eqnarray}\label{relacion entre n producto y (n,S)-producto}
   a_{\lambda}b&=&\hbox{Res}_z e^{\lambda z}Y(z)(a\otimes b)\nonumber\\
    &=& \hbox{Res}_z e^{(\lambda+T)z}Y^{op}(-z)(b\otimes a)\nonumber\\
     &=& -\hbox{Res}_z e^{-(\lambda+T)z}Y(z)S(z)(b\otimes a)\nonumber\\
     &:=& -b^{\mathcal{S}}_{-(\lambda+T)}a.
   \end{eqnarray}
   The collection of $n$-th product identities (\ref{axiom quantum n-product}) together (\ref{relacion entre n producto y (n,S)-producto}) are equivalent to:
   \begin{eqnarray}\label{desarollo jacobi}
      Y(a_{-\lambda}b, z)c&=& -Y(b^{\mathcal{S}}_{(\lambda-T)}a,z)c\nonumber\\
      &=&-\sum_{n \geq 0}Y(z)(b^{\mathcal{S}}_{(n)}a\otimes c)\frac{(\lambda-T)^{n}}{n!}\nonumber\\
      &=&-\sum_{n \geq 0}Y(z)^{\mathcal{S}}_{(n)}Y(z)(b\otimes a \otimes c)\frac{(\lambda-T)^{n}}{n!}\nonumber\\
      &=&-\sum_{n \geq 0}[\hbox{Res}_x \iota_{x,z}(x-z)^{n}Y(x)(1\otimes Y(z))(b\otimes a\otimes c)\nonumber\\
      &+& \iota_{z,x}(x-z)^{n}Y(z)(1\otimes Y(x))\mathcal{S}^{12}(z-x)(a\otimes b\otimes c)]\frac{(\lambda-T)^{n}}{n!}\nonumber\\
      &=&-\hbox{Res}_x e^{(\lambda-T)(x-z)}Y(x)(1\otimes Y(z))(b\otimes a \otimes c)\nonumber\\
      &+&\hbox{Res}_x e^{(\lambda-T)(x-z)}\sum_{i=1}^{r}e^{-x\partial_z}(f_{i}(z))Y(z)(1\otimes Y(x))(a^{i}\otimes b^{i} \otimes c)\nonumber\\
      &=&- e^{(-\lambda+T)z}[b_{(\lambda-T)}(Y(a, z) c)\nonumber\\
      &-&\sum_{i=1}^{r}\sum_{l\geq 0}(-\partial_z)^{(l)}(f_{i}(z))Y(a^{i},z)(b^{i}_{(\lambda-T, x^l)}c)].
      \end{eqnarray}
    Taking $\hbox{Res}_z e^{(\lambda+\mu)z},$ and changing $\lambda-T$ by $\alpha$ and $\mu+T$ by $\beta$, we obtain (\ref{jacoby identity (lambda,f) product}). Taking $\hbox{Res}_z z^{-1}$ in (\ref{desarollo jacobi}) and using $e^{(-\lambda+T)z}z^{-1}=z^{-1}+\int_{0}^{-\lambda+T}e^{\mu z}d\mu,$ we get
   \begin{eqnarray}
     (a_{-\lambda}b)_{\cdot}c &=& -b_{(\lambda-T)}(a_{\cdot}c)-\int_{0}^{T-\lambda}b_{(\lambda-T)}(a_{\mu}c)\,d\mu\nonumber\\
     &+&\sum_{i=1}^{r}\sum_{l\geq 0}(-1)^l[a^{i}_{(\cdot,(f_{i}(z))^{(l)})}(b^{i}_{(\lambda-T, x^l)}c)\nonumber\\
     &+&\int_{0}^{T-\lambda}  a^{i}_{(\mu,(f_{i}(z))^{(l)})}(b^{i}_{(\lambda-T, x^l)}c)\, d\mu].
     \end{eqnarray}
     This, together with (\ref{jacoby identity (lambda,f) product}), implies (\ref{relacion lambda y punto producto}) (after the substitution $\mu^{\prime}=\lambda+\mu-T$).

    Applying definitions of $(-1)$-st product and $Y^{op},$ due to  (\ref{skew-symmetry}) we get

   \begin{eqnarray}
   a_{\cdot}b&=&\hbox{Res}_z z^{-1}Y(z)(a\otimes b)\nonumber\\
    &=& \hbox{Res}_z z^{-1}e^{zT}Y^{op}(-z)(b\otimes a)\nonumber\\
     &=& \hbox{Res}_z z^{-1}e^{-zT}Y(z)S(z)(b\otimes a)\nonumber\\
     &=& \hbox{Res}_z z^{-1}Y(z)S(z)(b\otimes a)+\int_{0}^{-T}e^{\mu z}Y(z)S(z)(b\otimes a)d\mu\nonumber\\
     &=&b_{\cdot}^{\mathcal{S}}a+ \int_{0}^{-T}b_{\mu}^{\mathcal{S}}a\, d\mu.
   \end{eqnarray}
     Then this equation together the quantum  $(-1)$- product we get

   \begin{eqnarray}\label{relacion punto producto y S producto}
   Y(a\cdot b, z)c&=&Y(z)(b^{\mathcal{S}}_{(-1)}a\otimes c)+\int_{0}^{-T}Y(z)(b^{\mathcal{S}}_{\mu}a\otimes c)\, d\mu\nonumber\\
   &=&(Y(z)^{\mathcal{S}}_{(-1)}Y(z))(b\otimes a \otimes c)+\int_{0}^{-T}Y(z)(b^{\mathcal{S}}_{\mu}a\otimes c)\, d\mu\nonumber\\
   &=&Y(b,z)_{+}Y(a,z)c+\sum_{i=1}^{r}Y(a^{i},z)(b^{i}_{-\partial _z}c)(z^{-1}f_i(z))\nonumber\\
   &+&\int_{0}^{-T}Y(b^{\mathcal{S}}_{\mu}a,z)c\, d\mu.
   \end{eqnarray}
   Taking $\hbox{Res}_z e^{\lambda z}$  and using integration by parts, we get:
   \begin{eqnarray}
     (a_{\cdot}b)_{\lambda} &=& \hbox{Res}_z(e^{T\partial_{\lambda}}e^{\lambda z}b)_{\cdot}(a_{\lambda}c)+ \hbox{Res}_z\sum_{i=1}^{r}Y(a^i,z)(b^{i}_{\lambda-\partial_z}c (e^{\lambda z}z^{-1}f_i(z)))\nonumber\\
     &+&\int_{0}^{-T}(b^{\mathcal{S}}_{\mu}a)_{\lambda}c\, d\mu\nonumber\\
     &=& (e^{T\partial_{\lambda}}b)_{\cdot}(a_{\lambda}c)+ \hbox{Res}_z\sum_{i=1}^{r}Y(a^i,z)(b^{i}_{\lambda-\partial_z}c (z^{-1}f_i(z)
     +\int_{0}^{\lambda}f_i(z)e^{\mu z}\,d\mu)\nonumber\\
     &-&\int_{0}^{-T}(a_{-\mu-T}b)_{\lambda}c\, d\mu\nonumber\\
     &=& (e^{T\partial_{\lambda}}b)_{\cdot}(a_{\lambda}c) -\int_{0}^{-T}(a_{-\mu-T}b)_{\lambda}c\, d\mu
     + \sum_{i=1}^{r}[(e^{T\partial_{\lambda}}a^i)_{(\cdot, f_i(z))}(b^{i}_{\lambda}c)\nonumber\\
     &+&\sum_{l\geq 0}\int_{0}^{\lambda}a^i_{(\mu, (f_i(z))^{(l)})}(b^i_{(\lambda-\mu, x^l)}c\,)d\mu].
     \end{eqnarray}
  Due to (\ref{relacion punto producto y S producto}) and taking $\hbox{Res}_z z^{-1},$ we get
  \begin{eqnarray}
   (a_{\cdot} b)_{\cdot}c
   &=&\hbox{Res}_z z^{-1}Y(b,z)_{+}Y(a,z)_+c+ \hbox{Res}_z z^{-1}Y(b,z)_{+}Y(a,z)_-c\nonumber\\
    &+& \hbox{Res}_zz^{-1}\int_{0}^{-T}Y(b^{\mathcal{S}}_{\mu}a,z)c\, d\mu\nonumber\\
    &+& \sum_{i=1}^{r}\sum_{l\geq 0}(-1)^{l}\hbox{Res}_z z^{-1}Y(a^{i},z)(f_i(z))^{(l)}(\partial_z)^lY(b^{i}, z)_{-}c\nonumber\\
   &=&b_{\cdot}(a_{\cdot}c)+\hbox{Res}_z z^{-1})((e^{zT}-1)b)_{\cdot}(Y(a,z)_{-}c)
   +\int_{0}^{-T}(b_{\mu}^{\mathcal{S}}a)_{\cdot}c\, d\mu\nonumber\\
    &+&\sum_{i=1}^{r}\sum_{l\geq 0}\hbox{Res}_zz^{-1}((e^{zT}-1)a^{i})_{\cdot}((f_i(z))^{(l)}(\partial_z)^lY(b^{i}, z)_{-}c)\nonumber\\
    &+&\sum_{i=1}^{r}\sum_{l\geq 0}\hbox{Res}_zz^{-1}Y(a^{i}, z)_{-}((f_i(z))^{(l)}(\partial_z)^lY(b^{i}, z)_{-}c)\nonumber\\
    &=&b_{\cdot}(a_{\cdot}c)+\hbox{Res}_z\left(\int_{0}^{T}e^{\lambda z}d_{\lambda}b\right)_{\cdot}(Y(b,z)_{-}c)-\int_{0}^{-T}(a_{\mu-T}b)_{\cdot}c\,d\mu\nonumber\\
    &+&\sum_{i=1}^{r}\sum_{l\geq 0}\hbox{Res}_z\left(\int_{0}^{T}e^{\lambda z}d_{\lambda}a^{i}\right)_{\cdot}(f_i(z))^{(l)}(\partial_z)^lY(b^{i}, z)_{-}c)\nonumber\\
    &+&\sum_{i=1}^{r}\sum_{m,l\geq 0}(-1)^l a^{i}_{(\cdot, z^{m+1})}(b^{i}_{(.,\, D_{l}(f_i(z))z^{-m})}c),
   \end{eqnarray}
   which proves (\ref{asociatividad del producto punto}).
   \end{proof}

 Note that $V$ together with the $\lambda$-product is what \cite{BK} called \textit{conformal algebra}, and $V$ with $(\lambda, z^{-1})$-product is also a $\mathbb{K}[T]$-differential algebra with unit due to (\ref{identity}) and (\ref{derivation of the lambda productprima}).

 An important remark is that $V$ togheter with the $\lambda$-product is no longer a Leibnitz conformal algebra, since due to the braiding, the analog of the Jacoby identity involves $(\lambda, z^{-1})$-products.

With this in mind, we can prove the following
   \begin{theorem} Giving a braided state field correspondence $(V, |0\rangle, Y, \mathcal{S}),$ satisfying the $\mathcal{S}$-locality for $|0\rangle$ and the axiom of quantum $(n)$-product (\ref{quantum n-product}) implies to provide $V$ with a structure of a conformal algebra and  a structure of a $\mathbb{K}[T]$-differential algebra with a unit $|0\rangle,$ satisfying (\ref{jacoby identity (lambda,f) product})-(\ref{asociatividad del producto punto}).

   Conversely, given $V$ a topologically free $\mathbb{K}[[h]]$-module, a $\mathbb{K}[[h]]$- linear map $T$ and a braiding $\mathcal{S}.$  Assume that $V$  has a structure of conformal algebra and a structure of a $\mathbb{K}[T]$-differential algebra with a unit $|0\rangle,$ satisfying (\ref{jacoby identity (lambda,f) product})-(\ref{asociatividad del producto punto})
and $\mathcal{S}$ satisfies (\ref{Tx1 y S})- (\ref{1xT y S}), then $(V, |0\rangle, Y, \mathcal{S}),$ is a braided state field correspondence satisfying the axiom of quantum $(n)$-product, namely a strong braided field algebra.
   \end{theorem}

%XXXXXXXXXXXXX COMPLETAR PRUEBA
   \begin{proof}
   If $(V, |0\rangle, Y, \mathcal{S})$ is a braided state field algebra satisfying the axiom of quantum $n$-product, then by the above discussion we can define a $(\lambda, f)$-product on $V$ satisfying all the requirement.
  Conversely, given a $(\lambda, f)$-product we define a braided state field correspondence $Y$ by (\ref{definition state field correspondence}). In the proof of Lemma \ref{equations of a braided conformal algebra}, we have seen that the equations (\ref{jacoby identity (lambda,f) product})-(\ref{relacion lambda y punto producto}) are equivalent to the identities

  \begin{equation}\label{identity n-product}\hbox{Res}_z(Y({b_{n}^{\mathcal{S}}a},z)-Y(b,z)_{(n)}^{\mathcal{S}}Y(a,z))F=0, \qquad a, b \in V, n\geq 0,\, F=e^{\lambda z}\, \hbox{or}\, z^{-1}, \end{equation}
  while the equations (\ref{relation dot product and lambda product})-(\ref{asociatividad del producto punto}) are equivalent to the identities
   \begin{equation}\hbox{Res}_z\sum_{k\geq 0 }[(-\partial_z)^{(k)}Y(b_{(k-1)}^{\mathcal{S}}a,z)-\sum_{j=0}^{k}(-1)^{k}\partial_z^{(k-j)}Y(b,z)_{(k-1)}^{\mathcal{S}}\partial_z^{(j)}Y(a,z)]F=0,\end{equation}
   for $ a, b \in V$ and $ F=e^{\lambda z}\, \hbox{or}\, z^{-1}$.

   Due to Lemma \ref{T derivation of quantum n-product} and using translation invariance of $Y$, this identity is equivalent to
    \begin{equation}\label{identity (-1)-product}
    \hbox{Res}_z\sum_{k\geq 0 }[Y(b_{(k-1)}^{\mathcal{S}}a,z)-Y(b,z)_{(k-1)}^{\mathcal{S}}Y(a,z)](\partial_z)^{(k)}F=0,\end{equation} $\, a, b \in V,\, F=e^{\lambda z}\, \hbox{or}\, z^{-1} .$

    Using the translation invariance of $Y$ and integration by parts, we see that  identity (\ref{identity n-product}) holds also with $F$ replaced with $\partial_z F.$ Hence equations (\ref{identity n-product}) and (\ref{identity (-1)-product}) hold for all $F=z^{l}, l <0.$ For $F=e^{\lambda z},$ taking coefficients at power of $\lambda$ shows that they are satisfied also for $F=z^{l}, l\geq 0.$ This implies the $n$-th quantum product axioms for $n\geq -1.$ The the proof remains the same that proof of Theorem 4.4 \cite{BK}.
     \end{proof}

\section{Quantum conformal algebra}
%
%XXXXXXXXXXXXXXXXXXXXXXXXXXXXXXXXXXX

In this section, based on  what we have seen in Section 3, we aim to give a definition of braided conformal algebra. Until now, we didn't ask any furhter structure for the braiding $S$ besides (\ref{Tx1 y S}) and (\ref{1xT y S}). We have the following results.

\begin{proposition} If the hexagon relation
\begin{equation}\label{hexagon relation}
\mathcal{S}(x)(Y(z)\otimes 1)=(Y(z)\otimes 1)\mathcal{S}^{23}(x)\iota_{x,z}\mathcal{S}^{13}(x+z)
\end{equation}
holds in a braided state field correspondence, then we have that:
\begin{eqnarray}
\label{compatibility S with lambda, f product}
\mathcal{S}(x)(a_{(\lambda, f)}b\otimes c)&=&\sum_{l \geq 0}\partial_{\lambda}^l((\cdot_{(\lambda, f)}\cdot)\otimes 1)\mathcal{S}^{23}(x)\partial_{x}^{(l)}\mathcal{S}^{13}(x)(a\otimes b \otimes c)\\
&=&e^{\partial_{\lambda} \partial_{x_1}}((\cdot_{(\lambda, f)}\cdot)\otimes 1)\mathcal{S}^{23}(x)\mathcal{S}^{13}(x_1)(a\otimes b \otimes c)|_{x_1=x}
\end{eqnarray}
%\begin{equation}\label{Compatibility S with dot product}
%\mathcal{S}(x)(a_{(\cdot, f)}b\otimes c)=\sum_{l\geq 0}(\cdot_{(\cdot, f)}\cdot)\mathcal{S}^{23}(x)\partial_{x}^{(l)}\mathcal{S}^{13}(x)(a\otimes b \otimes c).
%\end{equation}
for all $a, b, c \in V$.% f \in \mathbb{K}((z))[[h]].$
\end{proposition}
\begin{proof}
Applying the definition of $(\lambda, f)$-product and using the hexagon relation (\ref{hexagon relation}), definition of S, Taylor expansion and change of variables we get,
 \begin{eqnarray}
%% \nonumber % Remove numbering (before each equation)
  \mathcal{S}(x)(a_{(\lambda, f)}b\otimes c) &=& \mathcal{S}(x)\hbox{Res}_z e^{\lambda z}f(z)(Y(z)\otimes 1)(a\otimes b \otimes c)\nonumber \\
 &=&\hbox{Res}_z e^{\lambda z}f(z)\mathcal{S}(x)(Y(z)\otimes 1)(a\otimes b \otimes c)\nonumber\\
  &=&\hbox{Res}_z e^{\lambda z}f(z)(Y(z)\otimes 1)\mathcal{S}^{23}(x)\iota_{x,z}\mathcal{S}^{13}(x+z)(a\otimes b \otimes c)\nonumber\\
&=&\sum_{i,j \in \mathbb{Z}}\hbox{Res}_z e^{\lambda z}f(z)(Y(z)\otimes 1)h_i(x)\iota_{x,z}g_j(x+z)(a^{(j)}\otimes b^{(i)}\otimes (c^{(j)})^{(i)})\nonumber\\
&=&\sum_{i,j \in  \mathbb{Z}}h_i(x)\hbox{Res}_z e^{\lambda z}f(z)(e^{z\partial_x }g_j(x))(Y(a^{(j)},z) b^{(i)}\otimes (c^{(j)})^{(i)})\nonumber\\
&=&\sum_{i,j, m, r \in \mathbb{Z}}\sum_{k, l \in \mathbb{Z}_{\geq 0} }h_i(x)\lambda^{(k)}f_r\,\hbox{Res}_z g^{(l)}_j(x)a^{(j)}_{(m)} b^{(i)}\otimes (c^{(j)})^{(i)}z^{-m-1+k+l+r}\nonumber\\
&=&\sum_{i,j, r \in \mathbb{Z}}\sum_{k, l \in \mathbb{Z}_{\geq 0} }h_i(x)\lambda^{(k)}f_r g^{(l)}_j(x)a^{(j)}_{(k+l+r)} b^{(i)}\otimes (c^{(j)})^{(i)}\nonumber\\
&=&\sum_{i,j, r \in \mathbb{Z}}\sum_{k\geq l \in \mathbb{Z}_{\geq 0} }h_i(x)\lambda^{(k-l)}f_r g^{(l)}_j(x)a^{(j)}_{(k+r)} b^{(i)}\otimes (c^{(j)})^{(i)}\nonumber\\
&=&\sum_{l\geq 0}\partial_{\lambda}^{l}(\cdot_{(\lambda, f)}\cdot)\mathcal{S}^{23}(x)\partial_x^{(l)}\mathcal{S}^{13}(x)(a\otimes b \otimes c),
  \end{eqnarray}

where $S^{23}(x)(a\otimes b\otimes c)=\sum_{i} h_i(x) a\otimes b^i\otimes c^i$ and $S^{13}(x)(a\otimes b\otimes c)=\sum_{j} g_j(x) a^j\otimes b\otimes c^j$. \end{proof}

Similarly, we have the following results.

\begin{proposition} \label{Associtivity proposition}If the associativity relation holds, namely,  there exists $N \in \mathbb{Z}_{\geq 0}$ such that
\begin{equation*}\label{associativity relation}
\iota_{z,w}(z+w)^NY(z+w)((1\otimes Y(w)))(a\otimes b \otimes c)
=(z+w)^NY(w)(Y(z)\otimes 1)(a\otimes b \otimes c)\,
 \end{equation*}
$\hbox{mod}\,h^M,$ for any $a, b, c \in  V$ and  $M \in \mathbb{Z}_{\geq 0}$ in a (braided) state field correspondence, then
\begin{equation}
\label{associativity with lambda, f product}
\partial_{\lambda}^N a_{\lambda}(b_{\mu}c)=\partial_{\lambda}^N (a_{\lambda}b)_{\lambda+\mu}c,
\end{equation}
%\begin{equation}\label{Compatibility S with dot product}
%\mathcal{S}(x)(a_{(\cdot, f)}b\otimes c)=\sum_{l\geq 0}(\cdot_{(\cdot, f)}\cdot)\mathcal{S}^{23}(x)\partial_{x}^{(l)}\mathcal{S}^{13}(x)(a\otimes b \otimes c).
%\end{equation}
$\hbox{mod}\,h^M$, for all $a, b, c \in V.$
\end{proposition}

%\end{document}
\begin{proof} Changing $z+w$ by $x$ in the associativity relation we have
\begin{equation}\label{associativity bis}
x^NY(x)(1\otimes Y(w))(a\otimes b \otimes c)=x^N Y(w)\iota_{x,w}(Y(x-w)\otimes 1)(a\otimes b \otimes c).
\end{equation}
Taking $\hbox{Res}_x\hbox{Res}_w e^{\lambda x}e^{\mu w}$ to the LHS of (\ref{associativity bis}) and using Lemma \ref{relation between lambda product y lambda f product} (a), we have
%\end{proof}
%\end{document}
\begin{equation}\label{RHS}\begin{aligned}
\hbox{Res}_x\hbox{Res}_w &e^{\lambda x}e^{\mu w}\, x^N Y(x)((1\otimes Y(w)))(a\otimes b \otimes c)\\
\,  &= (a_{(\lambda,x^N)}(b_{\mu} c))=\partial^N_{\lambda}(a_\lambda(b_\mu c)).
\end{aligned}
\end{equation}

Now, taking $\hbox{Res}_x\hbox{Res}_w e^{\lambda x}e^{\mu w}$ to the RHS of (\ref{associativity bis}),  using Taylor's formula (\ref{Taylor formula}), translation covariance and Lemma \ref{relation between lambda product y lambda f product} (a),

\begin{equation}\label{LHS}
\begin{aligned}
\hbox{Res}_x\hbox{Res}_w & e^{\lambda x}e^{\mu w}\,x^N Y(w)\iota_{x,w}(Y(x-w)\otimes 1)(a\otimes b \otimes c)\\
\, &=\hbox{Res}_x\hbox{Res}_w e^{\lambda x}e^{\mu w}\,x^N Y(w)e^{-w\partial_x}(Y(x)\otimes 1)(a\otimes b \otimes c)\\
\, &=\hbox{Res}_x\hbox{Res}_w e^{\lambda x}e^{\mu w}\,x^N Y(w)(Y(x)\otimes 1)((e^{-wT}a)\otimes b \otimes c)\\
\,&=\hbox{Res}_w e^{\mu w}Y( (e^{-wT}a)_{(\lambda,x^N)}b ,w)c\\
\,&=\hbox{Res}_w e^{\mu w} \partial_{\lambda}^N Y( (e^{-wT}a)_{\lambda}b ,w)c\\
\,&= \partial_{\lambda}^N\hbox{Res}_w e^{\mu w} e^{w\lambda}Y(a_{\lambda}b ,w)c\\
\,&=\partial_{\lambda}^N (a_{\lambda}b)_{\lambda+\mu}c.
\end{aligned}
\end{equation}

Equating (\ref{RHS}) and (\ref{LHS}), we finish the proof.\end{proof}

Now, we will show a similar resut but for the quasi-associativity (\ref{hexagonal relation}).

\begin{proposition} Let $V$ be a braided state field correspondance. Suppose that  for every $a,\,b,\,c\in V$ and $M\in\ZZ_{\geq 0}$ there exists $N\geq 0$ such that
$$\iota_{z,w}((z+w)^NY(z+w)(1\otimes Y(w))\mathcal{S}^{23}(w)\mathcal{S}^{13}(z+w)(a\otimes b\otimes c))$$
\begin{equation}\label{hexagonal relation}
=(z+w)^NY(w)\mathcal{S}(w)(Y(z)\otimes 1)(a\otimes b \otimes c)\quad \hbox{mod}\, h^{M}.
\end{equation}
holds mod$h^M$. Then
\begin{equation}\label{quasi assoc version lambda f}
\sum_{i,j}\partial_{\sigma}^Na^j_{(\sigma, g_j)} (b^i_{(-\lambda+\mu, h_i)} (c^j)^i)|_{\sigma=\lambda}=(\partial_\lambda+\partial_\mu)^N\sum_{r} (a_{\lambda}b)^r_{(\mu, f_r)} c^r
\quad mod\,h^M,
\end{equation}
where
$$\begin{aligned}
S(x)(a_\lambda b)\otimes c & =\sum_r f_r(x) (a_\lambda b)^r\otimes c^r),\\
S^{13}(x)(a\otimes b\otimes c)&=\sum_j g_j(x) (a^j\otimes b \otimes c^j),\\
S^{23}(x)(a^j\otimes b\otimes c^j)&=\sum_i h_i(x) (a^j\otimes b^i \otimes (c^j)^i).\\
\end{aligned}
$$
\end{proposition}

\begin{proof} Taking $\hbox{Res}_x\hbox{Res}_w e^{\lambda z}e^{\mu w}$ to the LHS of (\ref{hexagonal relation}), using Taylor's formula (\ref{Taylor formula}) and integration by parts, we have
\begin{equation}
\begin{aligned}
\hbox{Res}_z&\hbox{Res}_w e^{\lambda z}e^{\mu w}  \iota_{z,w}((z+w)^NY(z+w)(1\otimes Y(w))\mathcal{S}^{23}(w)\mathcal{S}^{13}(z+w)(a\otimes b\otimes c))=\\
\,&=\hbox{Res}_z\hbox{Res}_w e^{\lambda z}e^{\mu w} (z+w)^N e^{w\partial_z}( Y(z))(1\otimes Y(w))\mathcal{S}^{23}(w)e^{w\partial_z}(\mathcal{S}^{13}(z))(a\otimes b\otimes c))\\
\,&=\sum_{i,j} \hbox{Res}_z\hbox{Res}_w e^{-w\partial_z}(e^{\lambda z} (z+w)^N) e^{\mu w}Y(z))(1\otimes Y(w))g_j(z) h_i(w) (a^j\otimes b^i \otimes (c^j)^i).
\end{aligned}
\end{equation}

It is straightforward that $e^{-w\partial_z}(z+w)^N=z^N$ and $e^{-w\partial_z}e^{\lambda z}=e^{-w\lambda}e^{\lambda z}$, thus

\begin{equation}
\begin{aligned}
& \sum_{i,j} \hbox{Res}_z\hbox{Res}_w e^{-w\partial_z}\left(e^{\lambda z} (z+w)^N\right) e^{\mu w}Y(z)(1\otimes Y(w))g_j(z) h_i(w) (a^j\otimes b^i \otimes (c^j)^i)\\
\,&=\sum_{i,j} \hbox{Res}_z\hbox{Res}_w e^{-w\lambda}e^{\lambda z} (z)^N e^{\mu w}Y(z)(1\otimes Y(w))g_j(z) h_i(w) (a^j\otimes b^i \otimes (c^j)^i)\\
\,&=\sum_{i,j}a^j_{(\lambda, x^Ng_j)} (b^i_{(-\lambda+\mu, h_i)} (c^j)^i)\\
\,&=\sum_{i,j}\partial_{\sigma}^Na^j_{(\sigma,g_j)} (b^i_{(-\lambda+\mu, h_i)} (c^j)^i)|_{\sigma=\lambda}.
\end{aligned}
\end{equation}

In the last equality we used Lemma \ref{relation between lambda product y lambda f product}(a). Now, lets take residues in the RHS of (\ref{hexagonal relation}) and use Lemma \ref{relation between lambda product y lambda f product}(a)again. Thus

\begin{equation}
\begin{aligned}
& \hbox{Res}_z\hbox{Res}_w e^{\lambda z}e^{\mu w} (z+w)^NY(w)\mathcal{S}(w)(Y(z)\otimes 1)(a\otimes b \otimes c)\\
\,&=\sum_r\sum_{k=0}^N\binom{N}{k}\hbox{Res}_w e^{\mu w}f_r(w) w^{N-k} (a_{(\lambda, z^k)}b)^r\otimes c^r\\
\,&=\sum_r\sum_{k=0}^N\binom{N}{k} (a_{(\lambda, z^k)}b)^r_{(\mu,f_r(w) w^{N-k} )}\otimes c^r\\
\,&=(\partial_\lambda+\partial_\mu)^N\sum_{r} (a_{\lambda}b)^r_{(\mu, f_r)} c^r.
\end{aligned}
\end{equation}
Equating mod$h^M$, we have the desired result.\end{proof}

Finaly, let us translate the condition $Y(z)S(z)=Y^{op}(z)$ to the $(\lambda,f)$-product.

\begin{proposition} \label{skewsymmetry proposition}Suppose we have a state-field correspondance $V$ where $Y(z)S(z)=Y^{op}(z)$ holds. Then, for $a$ and $b\in V$,
\begin{equation}\label{skewsymmetry}
-b_{-\lambda-T}\,a=\sum_i a^i_{(\lambda,f_i)} b^i,
\end{equation}
where $\mathcal{S}(z)(a\otimes b)=\sum_i f_i(z) a^i\otimes b^i$

\end{proposition}
\begin{proof} We have that
\begin{equation}
\begin{aligned}
\hbox{Res}_z e^{\lambda z}Y(z)S(z)(a\otimes b) &=\sum_i \hbox{Res}_z e^{\lambda z}f_i(z) Y(a^i,z) b^i\\
\, &=\sum_i a^i_{(\lambda,f_i)} b^i.
\end{aligned}
\end{equation}
On the other hand, using $Y(z)S(z)=Y^{op}(z)$,
\begin{equation}
\begin{aligned}
\hbox{Res}_z e^{\lambda z}Y(z)S(z)(a\otimes b) &=\hbox{Res}_z e^{\lambda z}Y^{op}(z)(a\otimes b)\\
\,&=\hbox{Res}_z e^{\lambda z}e^{T z}Y(-z)(b\otimes a)\\
\,&=-\hbox{Res}_z e^{(-\lambda-T) z}Y(z)(b\otimes a)\\
\,&=-b_{-\lambda-T}\,a,
\end{aligned}
\end{equation}finishing the proof.
\end{proof}

A braided vertex algebra where  the associativity relation holds, is called \textit{quantum vertex algebra.} (Cf. Definition 3.12, \cite{DGK}). In the Characterization Theorem (cf. Theorem 5.13,\cite{DGK}) they  proved, among other equivalences, that a quantum vertex algebra is a braided state field correspondence such that the associativity relation and $YS=Y^{op}$ holds. We have shown in the discussion before Lemma \ref{Derivada de f en lambda prod}, combined with  the fact that all $(\lambda,f)$ products can be rewritten in terms of $\lambda$-products and $(\lambda,z^{-1})$-products, that having a braided state field correspondence is the same of having topologically free $\mathbb{K}[[h]]$-module $V,$ together with a
$\mathbb{K}[[h]]$-linear map $T: V\rightarrow V$, a distinguished vector $|0\rangle$, a braiding $\mathcal{S}$ on $V$ and  linear maps $(\lambda,f) : V\otimes V \rightarrow \mathbb{K}[\lambda][[h]], a\otimes b \rightarrow a_{(\lambda,f)}b$ for $f\in \mathbb{K}((Z))[[h]]$, such that

   \begin{equation}
    |0\rangle_{(\lambda,z^{-1})}a=a= a_{(\lambda,z^{-1})}|0\rangle,
   \end{equation}
   \begin{equation}\label{derivation of the lambda product}
    T(a_{(\lambda, f)}b)=T(a)_{(\lambda, f)}b+a_{(\lambda, f)}T(b)\end{equation} and
     \begin{equation}\label{derivation of the lambda product1}  T(a)_{(\lambda, f)}b=-\lambda a_{(\lambda, f)}b-a_{(\lambda, f')}b
   \end{equation}
   for all $a, b \in V$.  Combining this with Porposition \ref{Associtivity proposition} and Proposition \ref{skewsymmetry proposition} we have the following.

   \begin{theorem} Let $V$ be topologically free $\mathbb{K}[[h]]$-module, together with a
$\mathbb{K}[[h]]$-linear map $T: V\rightarrow V$, a distinguished vector $|0\rangle$, a braiding $\mathcal{S}$ on $V$. Define in $V$   linear maps $(\lambda,f) : V\otimes V \rightarrow \mathbb{K}[\lambda][[h]], a\otimes b \rightarrow a_{(\lambda,f)}b$ for $f\in \mathbb{K}((z))[[h]]$, such that the equation above hold. Let $Y$ be a  topological state-field correspondence.

   The following statements are equivalent:

   (i) $(V,\,T, |0\rangle, Y,\mathcal{S})$ is a quantum vertex algebra.

   (ii)$(V,T,|0\rangle, (\cdot_{(\lambda,f)}\cdot), \mathcal{S})$ satisfies the equations:

   \begin{equation}
    |0\rangle_{(\lambda,z^{-1})}a=a= a_{(\lambda,z^{-1})}|0\rangle,
   \end{equation}
   \begin{equation}\label{derivation of the lambda product}
    T(a_{(\lambda, f)}b)=T(a)_{(\lambda, f)}b+a_{(\lambda, f)}T(b)\end{equation} and
     \begin{equation}\label{derivation of the lambda product1}  T(a)_{(\lambda, f)}b=-\lambda a_{(\lambda, f)}b-a_{(\lambda, f')}b
   \end{equation}
   for all $a, b \in V$, and
   \begin{equation}
-b_{-\lambda-T}\,a=\sum_i a^i_{(\lambda,f_i)} b^i,
\end{equation}
where $\mathcal{S}(z)(a\otimes b)=\sum_i f_i(z) a^i\otimes b^i$, and there exists $N>>0$ such that

\begin{equation}
\label{associativity with lambda, f product}
\partial_{\lambda}^N a_{\lambda}(b_{\mu}c)=\partial_{\lambda}^N (a_{\lambda}b)_{\lambda+\mu}c,
\end{equation}
%\begin{equation}\label{Compatibility S with dot product}
%\mathcal{S}(x)(a_{(\cdot, f)}b\otimes c)=\sum_{l\geq 0}(\cdot_{(\cdot, f)}\cdot)\mathcal{S}^{23}(x)\partial_{x}^{(l)}\mathcal{S}^{13}(x)(a\otimes b \otimes c).
%\end{equation}
$\hbox{mod}\,h^M$, for all $a, b, c \in V.$
   \end{theorem}

If was proved in Proposition 3.13  in \cite{DGK} that if  a braided vertex algebra satisfies the hexagon relation then the associativity relation holds.

Assume that we have a braided vertex algebra $V$ and the hexagon relation holds, thus we have  a quantum vertex algebra. I we also ask in $V$ the condition
$$[T\otimes 1, \mathcal{S}(x)]=-\partial_x \mathcal{S}(x)\quad\hbox{and}\quad [1\otimes T, \mathcal{S}(x)]=\partial _x \mathcal{S}(x),$$ (which hold, for instance, in what \cite{EK5} called non-degenerate quantum vertex algebra), and  consider here  the $\lambda$-product  above, we showed that $(V, T,S)$ together with the $\lambda$-product is a conformal algebra ( in the sense of \cite{BK}),  sitting inside our quantum vertex algebra such that  (64) holds. All these, leads us to the following definition.

\,

\begin{definition}
  A \textit{quantum conformal algebra} is a topologically free $\mathbb{K}[[h]]$-module $V,$ together with a
$\mathbb{K}[[h]]$-linear map $T: V\rightarrow V$, a braiding $\mathcal{S}$ on $V$, namely a map

$$\begin{aligned} S:V\otimes V&\to \mathbb{C}((x))\otimes V\otimes V\\
a\otimes b &\mapsto \sum_{i=1}^rf_i(x) a^i\otimes b^i,
\end{aligned}$$
and a linear map $\lambda : V\otimes V \rightarrow \mathbb{K}[\lambda], a\otimes b \rightarrow a_{\lambda}b$ called $\lambda$-product, such that: ($a,b,c\in V$)
\begin{itemize}
  \item [(i)] $[T\otimes 1, \mathcal{S}(x)]=-\partial_x \mathcal{S}(x)$ (left shift condition);
  \item [(ii)] $[1\otimes T, \mathcal{S}(x)]=\partial _x \mathcal{S}(x)$ (right shift condition);
  \item [(iii)] $T(a_{\lambda}b)=(Ta)_\lambda b+ a_{\lambda} (Tb), \qquad (Ta)_{\lambda}b=-\lambda a_{\lambda}b$;
  \item[(iv)]$\label{compatibility S with lambda, f product}
\mathcal{S}(x)(a_{\lambda}b\otimes c)=e^{\partial_{\lambda} \partial_{x_1}}((\cdot_{\lambda}\cdot)\otimes 1)\mathcal{S}^{23}(x)\mathcal{S}^{13}(x_1)(a\otimes b \otimes c)|_{x_1=x}$, (hexagon relation).

\end{itemize}

Moreover if we ask $a_\lambda(b_\mu c)=(a_\lambda b)_{\lambda +\mu}c$, we call $V$ \textit{associative} quantum conformal algebra.

In particular, when  $S:V\otimes V\to \mathbb{C}[[x]]\otimes V\otimes V$ and  the Jacobi identity
\begin{equation}(a_{-\alpha-T}b)_{\alpha+\beta}c= -b_{\alpha}(a_{(\beta) }c)
      +\sum_{i=1}^{r}\sum_{l\geq 0} f^{(l)}_i(\partial_{\beta})(-\partial_{\alpha})^la^{i}_{\beta}(b^{i}_{\alpha}c),\end{equation}
  holds it is called\textit{ Leibnitz} quantum conformal algebra. Moreover,  if the Jacobi identity together with the skew-symmetry $-b_{-\lambda-T}a=\sum_{i=1}^r f_i(\partial_{\lambda}) (a^i_\lambda b^i)$ hold we say that $V$ is a \textit{Lie} quantum conformal algebra.

  \end{definition}

  Note that, due to (48) and (75), combined with Proposition 5(a),  the Jacobi identity  and the skew-symmetry  hold when we are in the case when the quantum vertex algebra has a braiding $S:V\otimes V\to \mathbb{C}[[x]]\otimes V\otimes V$, as in the following example.

\subsection{ Lie quantum conformal algebra  associated with associative algebras of Zamolodchikov-Faddeev type}.

In \cite{Li1} the quantum vertex algebras associated with associative
algebras of Zamolodchikov-Faddeev type is constructed as follows. Let $H$ be a vector space equipped with a bilinear form $\langle\cdot, \cdot \rangle$ and let $\mathcal{S}(x)$ be a linear map from $H\otimes H $ to $H\otimes H \otimes\mathbb{C}[[x]].$ Let $T(H \otimes \mathbb{C}[t, t^{-1}])^+$ be the suspace spanned by the vectors $$(a^{(1)}\otimes t^{n_1})\cdots (a^{(r)}\otimes t^{n_r}) $$ for $r\geq 1,$ $a^{(i)}\in H,$ $n_i \in \mathbb{Z}$ with $n_1+\cdots n_r \geq 0.$ Set
\begin{equation}
J=T(H \otimes \mathbb{C}[t, t^{-1}])T(H \otimes \mathbb{C}[t, t^{-1}])^+,
\end{equation}
a left ideal of $T(H \otimes \mathbb{C}[t, t^{-1}]).$ We then set
 \begin{equation}
 \tilde{V}(H, \mathcal{S})=T(H \otimes \mathbb{C}[t, t^{-1}])/J,
 \end{equation}
  a left $T(H \otimes \mathbb{C}[t, t^{-1}])$-module. Recall that we denote $a(n)=a\otimes t^n$ for $n\in\mathbb Z$. Clearly, $\tilde{V}(H, \mathcal{S})$ is cyclic on the vector $\tilde{1}=1+J$ where $a(n)\tilde{1}=0$ for $a \in H, \, n\geq 0. $ Furthermore, for $a \in H, \, w \in \tilde{V}(H, \mathcal{S}),$ we have
  \begin{equation}
  a(m)w=0 \quad \hbox{for}\, m \, \hbox{sufficiently large.}
  \end{equation}
  Now we define $V(H, \mathcal{S})$ to be the quotient $T(H \otimes \mathbb{C}[t, t^{-1}])$-module of $\tilde{V}(H, \mathcal{S})$ module the following relations:
  \begin{equation}
  a(x_1)b(x_2)w-\sum_{i=1}^{r}\iota_{x_2, x_1}(f_i(x_2-x_1)b^{(i)}(x_2)a^{(i)}(x_1)w=x_2^{-1}\delta \left(\frac{x_1}{x_2}\right)\langle a, b\rangle w
  \end{equation}
  for $a, b \in H, \, w \in \tilde{V}(H, \mathcal{S}),$ where $\mathcal{S}(b \otimes a )=\sum_{i=1}^{r}b^{(i)}\otimes a^{(i)}\otimes f_i(x).$ Denote by $\boldsymbol{1}$ the image of $\tilde{1}$ in $V(H, \mathcal{S}).$

  %Let $U$ a finite-dimensional vector space and an (invertible) element
%  \begin{equation}
%    R(x)=\sum_{n \geq 0}R_n x^n \in (\hbox{End}(U))[[x]]
%  \end{equation}
%  sucha that $R_iR_j=R_jR_i$ for $i,j \geq 0$ with $R_0(= R(0))$ invertible. Then we have $R(x_1)R(x_2)=R(x_2)R(x_1).$
%  Define $R^{*}(x) \in (\hbox{End}U)^*[[x]]$ by
%  \begin{equation}
%  \langle R^*(x)u^*,u  \rangle=\langle u^*, R^{-1}(x)u\rangle
%  \end{equation}
%  for $u \in U, u^* \in U^*,$ where $\langle \cdot, \cdot\rangle$ it the standard pairing between $U^{*}$ and $U.$ then
%  \begin{equation}
%  \langle R^*(x)u^*, R(x)u  \rangle=\langle u^*, u\rangle
%  \end{equation}
%  for $u \in U, u^* \in U^*.$
%  Set
% \begin{equation}
%%
%  H=U\oplus U^{*}
%  \end{equation}
%  and equip $H$ with the skew-symmetric bilinear form $\langle \cdot, \cdot\rangle$ defined by
%  \begin{equation}
%    \langle u+u^*, v+v^*\rangle= \langle u^*, v\rangle-\langle v^*, v\rangle
%  \end{equation}
%  for $u, v \in U, u^*, v^{*} \in U^*.$ Define $\mathcal{S}(x):h\otimes H \rightarrow H \otimes H \otimes \mathbb{C}[[x]]$ by
%  \begin{eqnarray}
%  % \nonumber % Remove numbering (before each equation)
%    \mathcal{S}(x)(u\times v) &=& R(-x)u\otimes R^{-1}(x)v, \\
%     \mathcal{S}(x)(u^{*}\times v^{*})&=& (R^{*})^{-1}(-x)u^{*}\otimes R^{*}(x)v^{*}, \\
%    \mathcal{S}(x)(v^{*}\times u)&=& R^{*}(-x)v^{*}\otimes R(x)u, \\
%    \mathcal{S}(x)(u\times v^{*})&=& R^{-1}(-x)u\otimes (R^{*})^{-1}(x)v^{*}.
%    \end{eqnarray}
%  for $u, v \in U, u^*, v^{*} \in U^*.$
    The following Theorem  in \cite{Li1}, shows  the existence of  a non-degenerate quantum vertex algebra structure on $V(H, \mathcal{S}):$
    \begin{theorem} Let $H$ be a finite-dimensional vector space over $\mathbb{C}$ eqquiped with a bilinear form $\langle \cdot, \cdot \rangle$ and let $\mathcal{S}: H\otimes H \rightarrow H \otimes H \otimes \mathbb{C}[[x]]$ be a bilinear map with $\mathcal{S}(0)=1.$ Assume that $V(H, \mathcal{S})$ is of PBW type in the sense that $S(H\otimes t^{-1}\mathbb{C}[t^{-1}])$-module gr$V(H, \mathcal{S})=\bigsqcup_{k \geq 0} V(H, \mathcal{S})[k]/V(H, \mathcal{S})[k-1]$ is a free module, where $V(H, \mathcal{S})[k]$ is the span of the vectors $a^{(1)}(-m_1)\cdots a^{(r)}(-m_r)\boldsymbol{1},$ for $r \geq 0,$ $a^{(1)}, \cdots, a^{(r)} \in H, \, m_i \geq 1 $ with $m_1+\cdots m_r \geq k.$ Then there exist a unique weak quantum structure on $V(H, \mathcal{S})$ with $\boldsymbol{1}$
    as the vacuum vector such that
  \begin{equation}
  Y(a(-1)\boldsymbol{1}, x)=a(x) \quad \hbox{for } a \in H
  \end{equation}
  Furhemore $V(H, \mathcal{S})$ is a nondegenerate quantum vertex algebra.
  \end{theorem}

Consider $V(H, \mathcal{S})$ as  in the theorem above as the base vector space together with the $\lambda$-product defined as $a_{\lambda}b=\Res_z e^{\lambda z} Y(a,z)b$. Note in particular, that given $a\in H$ and $w\in V(H,\mathcal S)$, due to (90) we have that

$$
(a(-1)\boldsymbol{1})_\lambda w=\sum_m (a(m)w)\, \lambda^m,
$$
where the sum is finite since (88) holds. As we saw in Section 3, this will be a quantum conformal algebra siting inside the quantum vertex algebras associated with associative
algebras of Zamolodchikov-Faddeev type since: the non-degenerancy of the quantum vertex algebra gives us conditions (i) and (ii), namely the left and right shift conditions in Definition 1. We showed in Section 3, that with this definition of $\lambda$-product in a general quantum vertex algebra (iii) always holds. As the hexagon relation holds in a quantum vertex algebra, we have (iv) by Proposition 5, and since the braiding $\mathcal{S}: H\otimes H \rightarrow H \otimes H \otimes \mathbb{C}[[x]]$ we also have Jacobi and skew-symmetry. Thus we have what we called a Lie quantum conformal algebra.

\,

In general, each time you consider a non-degenerate quantum vertex algebra with a braiding $\mathcal{S}: H\otimes H \rightarrow H \otimes H \otimes \mathbb{C}[[x]]$, you will have siting inside a Lie quantum conformal algebra.

\section{Author's contributions}
All authors contributed equally to this work.

\section{Data availability}
Data sharing is not applicable to this article as no new data were created or analyzed in this
study.

\end{document}